\documentclass[oneside,english]{amsart}
\usepackage[T1]{fontenc}
\usepackage[latin9]{inputenc}
\usepackage{mathrsfs}
\usepackage{amstext}
\usepackage{amsthm}
\usepackage{amssymb}
\usepackage[all]{xy}

\makeatletter
\numberwithin{equation}{section}
\numberwithin{figure}{section}

\PassOptionsToPackage{cmtip}{xy}
\PassOptionsToPackage{nooc}{xy}
\SilentMatrices
\makeatletter
\AtBeginDocument{%
  \def\xP@doSpecialRotate@@#1#2#3#4{}%
  \def\doSpecialRotate@@#1#2#3#4{}%
}
\makeatother

\makeatother

\theoremstyle{plain}
\newtheorem{thm}{\protect\theoremname}
\newtheorem{prop}[thm]{\protect\propositionname}
\theoremstyle{remark}
\newtheorem{rem}[thm]{\protect\remarkname}
\theoremstyle{definition}
\newtheorem{defn}[thm]{\protect\definitionname}
\theoremstyle{plain}
\newtheorem{lem}[thm]{\protect\lemmaname}
\newtheorem{cor}[thm]{\protect\corollaryname}
\theoremstyle{definition}
\newtheorem{example}[thm]{\protect\examplename}
\usepackage{babel}
\providecommand{\corollaryname}{Corollary}
\providecommand{\definitionname}{Definition}
\providecommand{\examplename}{Example}
\providecommand{\lemmaname}{Lemma}
\providecommand{\propositionname}{Proposition}
\providecommand{\remarkname}{Remark}
\providecommand{\theoremname}{Theorem}

\begin{document}
\title{Pointfree convergence approach and pointfree index theory}
\address{Department of Mathematics, Kean University, Jersey City campus, 2039
Kennedy Blvd, Jersey City, NJ 07305 }
\address{Department of Science and Math, Arkansas State University, Campus Querétaro, Querétaro, México} 
\email{eapzfc@gmail.com}
\email{fmynard@kean.edu}
\author{Emilio Angulo-Perkins and Frédéric Mynard}

\maketitle
\global\long\def\G{\mathcal{G}}%
 
\global\long\def\F{\mathcal{F}}%
 
\global\long\def\H{\mathcal{H}}%
\global\long\def\Z{\mathcal{Z}}%
 
\global\long\def\L{\mathcal{L}}%
\global\long\def\U{\mathcal{U}}%
\global\long\def\W{\mathcal{W}}%
 
\global\long\def\E{\mathcal{E}}%
\global\long\def\B{\mathcal{B}}%
 
\global\long\def\A{\mathcal{A}}%
\global\long\def\D{\mathcal{D}}%
\global\long\def\O{\mathcal{O}}%
 
\global\long\def\N{\mathcal{N}}%
 
\global\long\def\X{\mathcal{X}}%
 
\global\long\def\lm{\lim\nolimits}%
 
\global\long\def\then{\Longrightarrow}%
\global\long\def\BaseD{\operatorname{B}_{\mathbb{D}}}%

\global\long\def\V{\mathcal{V}}%
\global\long\def\C{\operatorname{C}}%
\global\long\def\R{\operatorname{R}}%
\global\long\def\adh{\operatorname{adh}}%
\global\long\def\Seq{\operatorname{Seq}}%
\global\long\def\intr{\operatorname{int}}%
\global\long\def\cl{\operatorname{cl}}%
\global\long\def\inh{\operatorname{inh}}%
\global\long\def\id{\operatorname{id}}%
\global\long\def\diam{\operatorname{diam}\ }%
\global\long\def\card{\operatorname{card}}%
\global\long\def\T{\operatorname{T}}%
\global\long\def\S{\operatorname{S}}%
\global\long\def\I{\operatorname{I}}%
\global\long\def\AdhD{\operatorname{A}_{\mathbb{D}}}%
\global\long\def\UD{\operatorname{U}_{\mathbb{D}}}%
\global\long\def\K{\operatorname{K}}%
\global\long\def\LD{\operatorname{L}_{\mathbb{D}}}%
\global\long\def\BRD{\operatorname{B}_{\mathcal{R}(\mathbb{D})}}%
\global\long\def\Dis{\operatorname{Dis}}%
\global\long\def\Cha{\operatorname{Cha}}%
\global\long\def\Reg{\operatorname{Reg}}%

\global\long\def\fix{\operatorname{fix}}%
\global\long\def\Epi{\operatorname{Epi}}%

\global\long\def\cont{\mathscr{C}}%
\global\long\def\conv{\mathsf{Conv}}%
\global\long\def\prtop{\mathsf{PrTop}}%
\global\long\def\Top{\mathsf{Top}}%
\global\long\def\pstop{\mathsf{PsTop}}%
\global\long\def\Cap{\mathsf{Cap}}%
\global\long\def\prap{\mathsf{PrAp}}%
\global\long\def\ap{\mathsf{Ap}}%
\global\long\def\psap{\mathsf{PsAp}}%
\global\long\def\Cconv{\mathsf{C^{Conv}}}%

\begin{abstract}
Given a category $\mathsf{C}$ of lattices, we introduce the category
$\mathsf{C^{Cap}}$, a pointfree version of the category $\Cap$ of
convergence approach spaces, that contains $\Cap$ coreflectively,
as well as the category $\mathsf{C^{Conv}}$, of pointfree convergence
spaces introduced by Goubault-Larrecq and Mynard, both reflectively
and coreflectively. We investigate pointfree counterparts of pseudo-approach,
pre-approach, and approach spaces along two natural axes---depth
and adherence-determinedness---which coincide in the spatial case
but diverge in general lattices, and we determine when the point--powerset
adjunction restricts to these subcategories. Unlike in $\mathsf{C^{Conv}}$,
the category $\ap$ of approach spaces finds a natural analog $\mathsf{C^{Ap}}$
in $\mathsf{C^{Cap}}$ and ``index theory'' translating traditional
topological or metric results into generalizing inequalities can be
conducted. We illustrate this with a result involving measure of continuity
and measure of compactness.
\end{abstract}

\section{Preliminaries}

\subsection{Set-theoretic and order-theoretic conventions}

If $X$ is a set, we will denote its powerset by $\mathbb{P}X$.

A partially ordered set $(L,\leq)$ is a lattice if all finite suprema
and infima (denoted with the symbols $\vee$ and $\wedge$ respectively)
exist. We will use $\bigvee$ and $\bigwedge$ for suprema and infima
over subsets of $L$, if they exist.

As usual, if $\alpha\in L$, $\uparrow\alpha=\{\beta\in L:\alpha\leq\beta\}$
and if $A\subset L$, $\uparrow A=\bigcup_{\alpha\in A}\uparrow\alpha$
and $\downarrow\alpha$ and $\downarrow A$ are defined dually. $A$
is called \emph{up-closed }(respectively, \emph{down-closed}) if $A=\uparrow A$
(respectively $A=\downarrow A$).We use interval notations when appropriate,
so that, for instance, given $\alpha,\beta\in L$, the notation $[\alpha,\beta]$
means $\uparrow\alpha\cap\downarrow\beta$. If $L$ has a least element
it is denoted by $\bot_{L}$ and if it has a greatest element, it
is denoted $\top_{L}$.

A non-empty subset $\F$ of $L$ is a \emph{filter }if it is up-closed
and closed under finite infima, that is, $a,b\in\F$ implies $a\wedge b\in\F$.
A filter $\F$ is \emph{proper }if it is a proper subset of $L$.
Note that given $\ell\in L$, $\uparrow\ell$ is a filter, called
\emph{principal filter of $\ell$. }Let $\mathbb{\bar{F}}L$ denote
the set of filters on $L$, while $\mathbb{F}L$ denotes the set of
\emph{proper }filters on $L$. Hence the set of traditional set-theoretic
filters on $X$ is $\mathbb{FP}X$, where $\mathbb{P}X$ is ordered
by inclusion. To make it easier to identify set-theoretic filters,
if $A\subset X$, we will write $A^{\uparrow}$ for $\uparrow A$
performed in $\mathbb{P}X$, that is, $\{B\subset X:A\subset B\}$.

We say that two subsets $A,B$ of a lattice $L$ \emph{mesh,} in symbols
$A\#B$, if $a\wedge b\neq\bot_{L}$ for every $a\in A$ and $b\in B$.
Moreover, we will use the notation
\[
A^{\#}=\{\ell\in L:\{\ell\}\#A\}.
\]
In particular, $\A,\B\subset\mathbb{P}X$ mesh if $A\cap B\neq\emptyset$
for every $A\in\A$ and $B\in\B$, and then $\A^{\#}=\{H\subset X:\forall A\in\A,A\cap H\neq\emptyset\}$. 

A filter $\F$ on a lattice $L$ is \emph{prime} if $a\in\F$ or $b\in\F$
whenever $a\vee b\in\F$ and \emph{completely prime} if there is $a\in A\cap\F$
whenever $\bigvee A\in\F$. On the other hand, an element $\ell$
of a lattice $L$ is \emph{meet-prime }if $a\wedge b\leq\ell$ implies
that either $a\leq\ell$ or $b\leq\ell$ and \emph{join-prime }if
$\ell\leq a\vee b$ implies that either $\ell\leq a$ or $\ell\leq b$.
Maximal proper filters are also called \emph{ultrafilters }and we
denote by $\mathbb{U}L$ the set of ultrafilters on $L$. In particular,
$\mathbb{UP}X$ denotes the set of set-theoretic ultrafilters on $X$.
If $\F\in\mathbb{F}L$, we write $\beta(\F)=\{\U\in\mathbb{U}L:\F\subset\U\}$.

\subsection{Convergence spaces}

A \emph{convergence }$\xi$ on a set $X$ is a relation between points
of $X$ and (set-theoretic) filters on $X$, denoted $x\in\lm_{\xi}\F$
if $(x,\F)\in\xi$ and interpreted as $x$ is a limit point for $\F$
in $\xi$, satisfying
\begin{equation}
\tag{centered}x\in\lm_{\xi}\{x\}^{\uparrow}\label{eq:centeredconv}
\end{equation}
for every $x\in X$ and 
\begin{equation}
\tag{monotone}\F\leq\G\then\lm_{\xi}\F\subset\lm_{\xi}\G,\label{eq:convmontone}
\end{equation}
for all $\F,\G\in\mathbb{FP}X$. The pair $(X,\xi)$ is then called
a \emph{convergence space.}

A map $f:X\to Y$ between two convergence spaces $(X,\xi)$ and $(Y,\tau)$
is \emph{continuous}, in symbols $f\in\cont(\xi,\tau)$, if $f(x)\in\lm_{\tau}f[\F]$
whenever $x\in\lm_{\xi}\F$, where $f[\F]=\{f(F):F\in\F\}^{\uparrow}$
is the image filter. Let $\conv$ denote the category of convergence
spaces and continuous maps. If $\xi,\sigma$ are two convergences
on $X$, we say that $\xi$ is \emph{finer than $\sigma$ }or that
$\sigma$ is \emph{coarser than $\xi$ }if the identity map $\id_{X}\in\cont(\xi,\sigma)$,
that is, if $\lim_{\xi}\F\subset\lim_{\sigma}\F$ for every filter
$\F$ on $X$.

If a convergence additionally satisfies
\begin{equation}
\tag{finite depth}\lm_{\xi}(\F\cap\G)=\lm_{\xi}\F\cap\lm_{\xi}\G\label{eq:fintedepth}
\end{equation}
 for every $\F,\G\in\mathbb{FP}X$, we say that $\xi$ \emph{has finite
depth} (\footnote{Many authors include (\ref{eq:fintedepth}) in the axioms of a convergence
space. Here we follow \cite{DM.book}.}). 

A convergence satisfying the stronger condition that 
\begin{equation}
\tag{pretopology}\lm_{\xi}(\bigcap_{\D\in\mathbb{D}}\D)=\bigcap_{\D\in\mathbb{D}}\lm_{\xi}\D\label{eq:pretopdeep}
\end{equation}
 for every $\mathbb{D}\subset\mathbb{FP}X$ is called a \emph{pretopology}
(and $(X,\xi)$ is a \emph{pretopological space}). 

A subset $A$ of a convergence space $(X,\xi)$ is \emph{$\xi$-open
}if 
\[
\lm_{\xi}\F\cap A\neq\emptyset\then A\in\F,
\]
and $\xi$-\emph{closed} if it is closed for limits, that is,
\[
A\in\F\then\lm_{\xi}\F\subset A,
\]
equivalently,
\[
A\in\F^{\#}\then\lm_{\xi}\F\subset A.
\]

Let $\O_{\xi}$ denote the set of open subsets of $(X,\xi)$ and let
$\O_{\xi}(x)=\{U\in\O_{\xi}:x\in U\}$. Similarly, let $\mathcal{C}_{\xi}$
denote the set of closed subsets of $(X,\xi)$. It turns out that
$\O_{\xi}$ is a topology on $X$\emph{. }Moreover,\emph{ }a topology
$\tau$ on a set $X$ determines a convergence $\xi_{\tau}$ on $X$
by
\[
x\in\lm_{\xi_{\tau}}\F\iff\F\geq\N_{\tau}(x),
\]
where $\N_{\tau}(x)$ denotes the neighborhood filter of $x$ for
$\tau$. In turn, $\xi_{\tau}$ completely determines $\tau$ because
$\tau=\O_{\xi_{\tau}}$, so that we do not distinguish between $\tau$
and $\xi_{\tau}$ and identify topologies with special convergences.
Moreover, a convergence $\xi$ on $X$ determines the topology $\O_{\xi}$
on $X$ which turns out to be the finest among the topologies on $X$
that are coarser than $\xi$. We call it the \emph{topological modification
of $\xi$ }and denote it $\T\xi$. Hence the category $\Top$ of topological
spaces and continuous maps is a concretely reflective subcategory
of $\conv$. Let $\cl_{\xi}$ denote the closure operator in $\T\xi$.

A convergence $\xi$ determines an \emph{adherence} operator on filters,
namely,
\begin{equation}
\adh_{\xi}\F=\bigcup_{\G\#\F}\lm_{\xi}\G=\bigcup_{\G\geq\F}\lm_{\xi}\G=\bigcup_{\U\in\beta(\F)}\lm_{\xi}\U.\label{eq:adherence}
\end{equation}

A convergence $\xi$ on $X$ is a \emph{pseudotopology }if 
\begin{equation}
\tag{pseudotopology}\lm_{\xi}\F=\bigcap_{\U\in\beta(\F)}\lm_{\xi}\U\label{eq:pseudotop}
\end{equation}
 for every filter $\F\in\mathbb{FP}X$. The full subcategory of $\conv$
formed by pseudotopological spaces (and continuous maps), denoted
$\pstop$, is concretely reflective and the reflector $\S$ can be
described on objects by
\begin{equation}
\lm_{\S\xi}\F=\bigcap_{\U\in\beta(\F)}\lm_{\xi}\U=\bigcap_{\G\#\F}\adh_{\xi}\G,\label{eq:Sreflector}
\end{equation}
so that adherences determine the pseudotopological reflection.

When restricted to principal filters, (\ref{eq:adherence}) defines
a \emph{principal adherence operator }$\adh_{\xi}:\mathbb{P}X\to\mathbb{P}X$
given by
\[
\adh_{\xi}A=\adh_{\xi}\{A\}^{\uparrow}=\bigcup_{A\in\G^{\#}}\lm_{\xi}\G=\bigcup_{A\in\G}\lm_{\xi}\G=\bigcup_{\U\in\beta A}\lm_{\xi}\U.
\]

The full subcategory of $\conv$ formed by pretopological spaces (and
continuous maps), denoted $\prtop$, is concretely reflective and
the reflector $\S_{0}$ can be described on objects by
\begin{equation}
\lm_{\S_{0}\xi}\F=\bigcap_{A\in\F^{\#}}\adh_{\xi}A,\label{eq:S0adh}
\end{equation}
so that the principal adherence determines the pretopological reflection.
In general,
\[
\adh_{\xi}A\subset\cl_{\xi}A
\]
but not conversely. In contrast to $\cl_{\xi}$, the principal adherence
is in general non-idempotent because $\adh_{\xi}A$ need not be closed.
A pretopology is a topology if and only if its principal adherence
operator is idempotent, in which case $\adh_{\xi}=\cl_{\xi}$. Moreover,
the reflector $\T$ from $\conv$ onto $\Top$ is given by
\begin{equation}
\lm_{\T\xi}\F=\bigcap_{A\in\F^{\#}}\cl_{\xi}A.\label{eq:Tclosure}
\end{equation}

\subsubsection{\label{subsec:Compactness-for-filters}Compactness for filters in
$\mathsf{Conv}$}

As we will discuss compactness in Section \ref{sec:Measure-of-compactness},
let us note that a convergence space is \emph{compact }if every ultrafilter
converges, equivalently if the adherence of every filter is non-empty.
Generalizations in several directions have proved to be unifying notions
in the context of convergence spaces. Given a class $\mathbb{D}$
of filters (\footnote{that is, assignment $\mathbb{D}X\subset\mathbb{FP}X$ for every set
$X$, which to simplify we're going to assume only depends on the
set $X$ and not the convergence, though this can be usefully considered.}), a filter $\F$ on a convergence space $(X,\xi)$ and $\A\subset\mathbb{P}X$,
we say that $\F$ is $\mathbb{D}$-\emph{compact at $\A$ }if
\[
\D\#\F\then\adh_{\xi}\D\#\A,
\]
for every $\D\in\mathbb{D}X$ and \emph{nearly }$\mathbb{D}$-\emph{compact
at $\A$} if
\[
\D\geq\F\then\adh_{\xi}\D\#\A,
\]
for every $\D\in\mathbb{D}X$. To say that $F\subset X$ is compact
(resp. countably compact, Lindelöf) is to say that the principal filter
$\{F\}^{\uparrow}$ is $\mathbb{D}$-compact at itself, that is, at
$\{F\}$ (equivalently, nearly $\mathbb{D}$-compact at $\{F\}$),
where $\mathbb{D}$ is the class $\mathbb{F}$ of all filters (resp.,
the class $\mathbb{F}_{1}$ of countably based filters, the class
$\mathbb{F}_{\wedge1}$ of countably complete filters). It is relatively
compact if $\{F\}^{\uparrow}$ is compact at $\{X\}$. Accordingly,
we omit the prefix $\mathbb{F}$- when $\mathbb{D}=\mathbb{F}$. Note
that $x\in\lm_{\S\xi}\F$ if and only if $\F$ is compact at $\{x\}$
and $x\in\lm_{\S_{0}\xi}\F$ if and only if $\F$ is $\mathbb{F}_{0}$-compact
at $\{x\}$. See \cite{D.comp,myn.applofcompact,myn.relations} for
an account of the unifying power of these notions. 

Let's say that a class $\mathbb{D}$ of filters is \emph{transferable
}if $\D\in\mathbb{D}Y$ implies $f^{-}[\D]\in\mathbb{D}X$ for every
$f:X\to Y$. Note that the following is a simple generalization of
the fact that the continuous image of a compact set is compact (see,
e.g., \cite[Proposition IX.8.2, Corollary XV.3.3]{DM.book}):
\begin{prop}
\label{prop:contimageofcompact} If $\mathbb{D}$ is a transferable
class of filters, $(X,\xi)$ and $(Y,\tau)$ are convergence spaces,
$\F$ is a $\mathbb{D}$-compact filter at $A\subset X$ and $f\in\cont(\xi,\tau)$
then $f[\F]$ is $\mathbb{D}$-compact at $f(A)$.
\end{prop}

\begin{rem}
In particular, applied to principal filters, this yields the classical
fact that compactness, countable compactness and the Lindelöf property
are preserved under continuous maps. Applied to filters and $A=\{x\}$,
this yields that $\S$, $\S_{0}$ and other similarly defined maps
preserve continuity, that is, are functors. Though both sets of results
are easy to obtain, it is interesting to see them as two sides of
the same coin, under the umbrella of Proposition \ref{prop:contimageofcompact}.
\end{rem}

We refer the reader to \cite{DM.book} for a systematic study of convergence
spaces.

\subsection{Pointfree convergence spaces\label{subsec:Pointfree-convergence-spaces}}

Note that (\ref{eq:convmontone}) states that $\lim_{\xi}:\mathbb{FP}X\to\mathbb{P}X$
is monotone, which allows for a pointfree abstraction of $\conv$,
after observing that given two convergence spaces $(X,\xi)$ and $(Y\,\tau)$,
the condition $f\in\cont(\xi,\tau)$ rephrases as$f(\lm_{\xi}\F)\subset\lm_{\tau}f[\F]$,
equivalently,
\[
\lm_{\xi}\F\subset\varphi(\lm_{\tau}\varphi^{-1}[\F]),
\]
where $\varphi=f^{-1}:\mathbb{P}Y\to\mathbb{P}X$ is a morphism of
lattices.
\begin{defn}
\cite{FredetJean} Given a category $\mathsf{C}$ of lattices, a\emph{
convergence $\mathsf{C}$-object $(L,\lim)$ }is a $\mathsf{C}$-object
$L$ together with a monotone map $\lim:\mathbb{F}L\to L$. The objects
of the category $\Cconv$ are the convergence $\mathsf{C}$-objects
and the morphisms $\varphi:L\to L'$ are the $\mathsf{C}$-morphisms
that are \emph{continuous }in the sense that for every $\F\in\mathbb{F}L'$,
\begin{equation}
\lm_{L'}\F\leq\varphi(\lm_{L}\varphi^{-1}(\F)),\tag{ptfree continuity}\label{eq:ptfreeCont}
\end{equation}
 where $\varphi^{-1}(\F)=\{\ell\in L:\varphi(\ell)\in\F\}$.
\global\long\def\pt{\operatorname{pt}}%

If $\mathsf{C}$ denotes a category $\mathsf{C^{op}}$ denotes its
opposite category. The category $\conv$ embeds coreflectively into
$(\Cconv)^{op}$ when $\mathsf{C}$ is the category of frames or of
coframes: the \emph{powerset functor} $\mathbb{P}:\conv\to(\Cconv)^{op}$
sending $(X,\xi)$ to $(\mathbb{P}X,\lm_{\xi})$ and $f:(X,\xi)\to(Y,\tau)$
to $\text{\ensuremath{\mathbb{P}f=}}f^{-1}:\mathbb{P}Y\to\mathbb{P}X$
(in $\Cconv$) is then right-adjoint to the \emph{point-functor} $\pt:(\Cconv)^{op}\to\conv$
(the coreflector) where the underlying set of $\pt L$, the set of
``points'' of $L$, is the set of $\Cconv$-morphisms from $L$
to $\mathbb{P}(\{*\})$, hence depends on the choice of $\mathsf{C}$.
The convergence structure on $\pt L$ is given by 
\begin{equation}
\tag{Conv on pt}\lm_{\pt L}\F=\left(\lm_{L}\F^{\circ}\right)^{\bullet},\label{eq:ptconv}
\end{equation}
where $\ell^{\bullet}=\{\varphi\in\pt L:\varphi(\ell)=\{*\}\}$ and
$\F^{\circ}=\left\{ \ell\in L:\ell^{\bullet}\in\F\right\} $. Finally,
if $\varphi\in\Cconv(L,L')$ then $\pt\varphi:\pt L'\to\pt L$ is
defined by $\pt(\varphi)(f)=f\circ\varphi$.
\end{defn}

\begin{rem}
\label{rem:ptPXinCconv} Identifying a point $\varphi\in\pt L$ with
the filter $\F=\varphi^{-1}(\{*\})\in\mathbb{F}L$, points can be
seen as special filters on $L$. The continuity condition (\ref{eq:ptfreeCont})
on $\varphi$ yields $\lim\F\in\F$, while the fact that $\varphi$
is a $\mathsf{C}$-morphism adds another condition. In particular,
if $\mathsf{C}$-morphisms preserves arbitrary suprema, for instance
in the category of frames, then $\F$ is a completely prime filter.
If $\mathsf{C}$-morphisms preserves arbitrary infima, for instance
in the category of coframes, then $\F$ is the principal filter of
a join-prime element. Note that for these two cases, the only such
filters in $\mathbb{P}X$ are principal ultrafilters of singletons,
so that the underlying set of $\pt\mathbb{P}X$ can be identified
with $X$ for such base categories $\mathsf{C}$, since $x\in\lim\{x\}^{\uparrow}$
for every $x\in X$. See \cite{FredetJean} and \cite{MYNARD2022108226}
for details.
\end{rem}

\subsection{\label{subsec:Complete-Heyting-algebras}Complete Heyting algebras
and grills}

Let us focus on the case where $\mathsf{C}$ is the category $\mathsf{CH}$
of complete Heyting algebras. The objects of $\mathsf{Heyt}$ are
frames, i.e., complete lattices $L$ with
\begin{equation}
a\wedge\bigvee_{b\in B}b=\bigvee_{b\in B}(a\wedge b)\label{eq:frame}
\end{equation}
 for every $a\in L$ and $B\subset L$. As $a\wedge\cdot$ preserves
all suprema, it has a right adjoint $a\to\cdot$, so that
\begin{equation}
a\wedge x\leq y\iff x\leq a\to y\label{eq:implication}
\end{equation}
 for every $a,x,y\in L$. The element $a^{*}:=a\to0$ of $L$ is called
the \emph{pseudocomplement of $a$. }The morphisms of $\mathsf{CH}$
preserve arbitrary joins and meets, and the implication $\to$. Note
that a morphism $\varphi\in\mathsf{CH}(L,L')$ preserves all infima
(in particular, $\varphi(1_{L})=1_{L'}$) and therefore admits a left-adjoint
$\varphi_{\ell}:L'\to L$, and preserves all suprema (in particular
$\varphi(0_{L})=0_{L'}$) and therefore admits a right-adjoint $\varphi_{r}:L'\to L$
so that
\begin{equation}
\varphi_{\ell}(m')\leq m\iff m'\leq\varphi(m),\label{eq:leftadjoint}
\end{equation}
and
\begin{equation}
\varphi(m)\leq m'\iff m\leq\varphi_{r}(m'),\label{eq:rightadjoint}
\end{equation}
whenever $m\in L$ and $m'\in L'$.

In this context, we say that two subsets $A$ and $B$ of a complete
lattice $L$ \emph{mesh, }in symbols $A\#B$, if $a\wedge b\neq0_{L}$
whenever $a\in A$ and $b\in B$. We will also consider the \emph{grill
of $A\subset L$ }given by $A^{\#}=\{\ell\in L:\{\ell\}\#A\}.$ As
in a frame $a\wedge b=0$ if and only if $b\leq a^{*}$, we have that
\begin{equation}
\ell\in A^{\#}\iff\ell^{*}\notin\uparrow A\label{eq:grillpseudoc}
\end{equation}
 for every $A\subset L$ if $L$ is a frame. 
\begin{lem}
\label{lem:grillofunion} If $L$ is a distributive lattice, $\H\in\mathbb{F}L$
and $A,B$ are two upsets in $L$ then 
\[
\H\#(A\cap B)\iff\H\#A\text{ or }\H\#B.
\]
\end{lem}

\begin{proof}
If $A\not\subset\H^{\#}$ and $B\not\subset\H^{\#}$ there is $a\in A$
and $b\in B$ with $a\notin\H^{\#}$ and $b\notin\H^{\#}$, that is,
there is $h_{1},h_{2}\in\H$ with $a\wedge h_{1}=b\wedge h_{2}=0$.
Moreover, $a\vee b\in A\cap B$ because $A$ and $B$ are upsets,
and $h=h_{1}\wedge h_{2}\in\H$ because $\H$ is a filter. As 
\[
h\wedge(a\vee b)=(h\wedge a)\vee(h\wedge b)\leq(h_{1}\wedge a)\vee(h_{2}\wedge b)=0,
\]
we conclude that $A\cap B\not\subset\H^{\#}$. The reverse implication
is clear.
\end{proof}

Moreover, \cite[Lemma 6.19]{FredetJean} states (in slightly different
terms, but equivalently) that if $\varphi:L\to L'$ is a morphism
of coframes, $M\subset L$ and $N\subset L'$ then
\begin{equation}
\varphi(M)\#N\then M\#\varphi_{\ell}(N).\label{eq:grillcoframe}
\end{equation}
It turns out that the converse is true for morphisms of complete Heyting
algebras. Indeed, in the category $\mathsf{CH}$, morphisms preserves
implications, so that 
\begin{equation}
\varphi\in\mathsf{CH}(L,L')\then\varphi(\ell^{*})=(\varphi(\ell))^{*},\label{eq:CHpseudocompl}
\end{equation}
for every $\ell\in L$. As a result:
\begin{lem}
\label{lem:Heygrill} Let $\varphi\in\mathsf{CH}(L,L')$, $M\subset L$
and $N\subset L'$. Then 
\begin{equation}
\varphi(M)\#N\iff M\#\varphi_{\ell}(N).\label{eq:grillCH}
\end{equation}
\end{lem}

\begin{proof}
Let $m\in M$ and $n\in N$. Then
\[
\varphi(m)\wedge n=0_{L'}\iff n\leq(\varphi(m))^{*}=\varphi(m^{*})\iff\varphi_{\ell}(n)\leq m^{*}\iff\varphi_{\ell}(n)\wedge m=0_{L}.
\]
\end{proof}

Consider now $\mathsf{CH^{Conv}}$. Note that $\mathbb{P}:\mathsf{Conv}\to\mathsf{(CH^{Conv})^{op}}$
is still a functor as $\mathbb{P}f=f^{-1}:\mathbb{P}Y\to\mathbb{P}X$
preserves arbitrary unions, intersections, and complements, hence
the implication $A\to B=(X\setminus A)\cup B$. On the other hand,
if $L\in\mathsf{Ob(CH^{Conv})}$ then $\pt L$ can be identified with
filters on $L$ that are principal filters of completely join-prime
(or compact) elements such that either $a$ or $a^{*}$ is in the
filter, whenever $a\in L$. In particular, such filters in $\mathbb{P}X$
are the principal filters of singletons, so that $\pt\mathbb{P}X$
can be identified with $X$.

\subsection{The quantale $V$\label{subsec:V}}

$V$ will denote a complete and completely distributive lattice with
top element $1=\top_{V}$ and bottom element $0=\bot_{V}$, equipped
with a quantale operation $\otimes:V\times V\to V$ which is commutative
and unital, with unit $1$ (an integral quantale). Hence,
\begin{equation}
a\otimes\bigvee_{b\in B}b=\bigvee_{b\in B}a\otimes b\label{eq:quantalesup}
\end{equation}
for every $a\in V$ and $B\subset V$, and 
\begin{equation}
\forall a\in V,\;a\otimes1=a.\label{eq:atimes1}
\end{equation}

We also require that $1$ does not have a unique immediate predecessor,
that is,
\begin{equation}
1=\bigvee_{\epsilon\in V\setminus\{1\}}\epsilon.\label{eq:density}
\end{equation}

Standard examples include $V=[0,1]$ with $\otimes=\times$ or $\otimes=\wedge$.
\begin{rem}
\label{rem:quantaleiso} Note that $-\ln:([0,1],\times)\to([0,\infty],+)^{op}$
is a quantale isomorphism, and thus the classical setting of approach
spaces is covered, albeit with the opposite order on $[0,\infty]$. 

We gather below straightforward facts about $(V,\otimes)$:
\end{rem}

\begin{lem}
\label{lem:basicotimes} Let $(V,\otimes)$ be as above. 
\begin{enumerate}
\item \label{ittimes1:below}$a\otimes v\leq a\wedge v$ for every $a,v\in V$;
\item \label{ittimes0}$a\otimes0=0$ for every $a\in V$;
\item \label{ittimes:orderpreserving}If $a,b,c,d\in V$ then
\[
a\leq b\text{ and }c\leq d\then a\otimes c\leq b\otimes d.
\]
\item \label{itsupoftimes}if $v\in V$ then
\[
v=\bigvee_{\epsilon\in V\setminus\{1\}}v\otimes\epsilon.
\]
\end{enumerate}
\end{lem}

In $(V,\otimes)$, for each $v\in V$, the map $v\otimes-:V\to V$
preserves all suprema, hence has a right adjoint $-\oslash v:V\to V$
defined by 
\begin{equation}
v\otimes x\leq y\iff x\leq y\oslash v,\label{eq:adjointtensor}
\end{equation}
which preserves all infima, namely, 
\begin{equation}
\left(\bigwedge_{b\in B}b\right)\oslash v=\bigwedge_{b\in B}(b\oslash v).\label{eq:oslahinf}
\end{equation}

Note as a result the following straightforward properties of $\oslash$:
\begin{lem}
\label{lem:basicoslah} For every $v,x,t\in V$,
\begin{enumerate}
\item \label{it1:tidvgreaterthant}$t\oslash v\geq t$;
\item \label{it2:divis1}$v\oslash t=1$ if and only if $v\geq t$ ( in
particular, $1\oslash v=1$, $v\oslash v=1$, and $v\oslash0=1$.)
\item \label{it3:divby1}$v\oslash1=v$;
\item \label{it4:orderreversing}
\[
x\leq t\then v\oslash x\geq v\oslash t;
\]
\item \label{it5:divbysup}
\[
v\oslash\bigvee_{a\in A}a=\bigwedge_{a\in A}(v\oslash a);
\]
\item \label{it6:divofdiv}
\[
(v\oslash t)\oslash(v\oslash x)\geq x\oslash t;
\]
\item \label{it7:infofdiv}
\[
v=\bigwedge_{\epsilon\in V\setminus\{1\}}v\oslash\epsilon.
\]
\end{enumerate}
\end{lem}

\section{The category $\mathbf{\mathsf{C}^{\protect\Cap}}$}

\subsection{Reframing $\protect\Cap$}

Convergence approach spaces were introduced in \cite{Lowen88}. Traditionally,
a \emph{convergence approach space }$(X,\lambda)$ is given by a function
$\lambda:\mathbb{FP}X\to[0,\infty]^{X}$ satisfying 
\[
\F\leq\G\then\lambda(\F)\geq\lambda(\G)
\]
for the pointwise order of $[0,\infty]^{X}$ induced by the usual
order on $[0,\infty]$, and 
\[
\lambda(\{x\}^{\uparrow})(x)=0
\]
and a morphism (then called \emph{contraction}) is a function $f:(X,\lambda_{X})\to(Y,\lambda_{Y})$
satisfying 
\[
\lambda_{Y}(f[\F])\circ f\leq\lambda_{X}(\F)
\]
 for this pointwise order, for every $\F\in\mathbb{FP}X$.
\begin{rem}
In fact, most authors additionally require that 
\begin{equation}
\lambda(\F\cap\G)=\lambda(\F)\vee\lambda(\G)\label{eq:finitedeptholdcap}
\end{equation}
for every $\F$ and $\G$ in $\mathbb{FP}X$, a property that is inessential
for most purposes and that we therefore will not require and call
\emph{finite depth}.
\end{rem}

In view of Remark \ref{rem:quantaleiso}, we prefer to reframe the
structure as a monotone map $\lambda:\mathbb{FP}X\to[0,1]^{X}$ for
the pointwise order, with 
\begin{equation}
\lambda(\{x\}^{\uparrow})(x)=1\label{eq:CAPcentered}
\end{equation}
 for every $x\in X$, and morphisms $f:(X,\lambda_{X})\to(Y,\lambda_{Y})$
satisfying 
\[
\lambda_{Y}(f[\F])\circ f\geq\lambda_{X}(\F)
\]
for every $\F\in\mathbb{FP}X$. 
\global\long\def\cont{\mathscr{C}}%
 Let $\cont(X,Y)$ denote the set of all morphisms. The resulting
category is denoted $\Cap$.

Note that $\Cap$ is a topological category, where initial structures
are as follows: If $f_{i}:X\to(Y_{i},\lambda_{i})$ then the $\Cap$-initial
structure on $X$ is given by 
\begin{equation}
\lambda_{X}(\F)(x)=\bigwedge_{i\in I}\lambda_{i}(f_{i}[\F])(f_{i}(x)).\label{eq:initialCAP}
\end{equation}

Analogously to (\ref{eq:pretopdeep}), a convergence approach space
$(X,\lambda)$ is a called a \emph{pre-approach space }(e.g., \cite{Lowen97})
if for any $\mathbb{D}\subset\mathbb{FP}X$,
\begin{equation}
\lambda(\bigcap_{\D\in\mathbb{D}}\D)(\cdot)=\bigwedge_{\D\in\mathbb{D}}\lambda\D(\cdot).\tag{PrAp}\label{eq:preAP}
\end{equation}

If (\ref{eq:preAP}) is only true for subsets $\mathbb{D}$ of $\mathbb{UP}X$
of the form $\beta(\F)$, that is, if 
\begin{equation}
\tag{PsAp}\lambda(\F)(\cdot)=\bigwedge_{\U\in\beta(\F)}\lambda\U(\cdot)\label{eq:psap}
\end{equation}
for every $\F\in\mathbb{F}X$, then $(X,\lambda)$ is a called a \emph{pseudo-approach
space }(e.g.,\cite{lowe88}), generalizing pseudotopologies (compare
(\ref{eq:pseudotop})). We denote the corresponding full reflective
subcategories of $\Cap$ by $\mathbf{\prap}$ and $\psap$ respectively.
Though classical approach theory as in \cite{AP.book,indextheory}
is concerned with the subcategory $\mathsf{Ap}$ of $\prap$ formed
by \emph{approach spaces}, we will only consider them later after
we introduce appropriate notions of adherence.

Note that $\conv$ embeds both reflectively and coreflectively in
$\Cap$: A $\conv$-object $\left(X,\lim\right)$ uniquely determines
a $\Cap$-object $e(X)=\left(X,\lambda_{\lm}\right)$ defined by 
\begin{equation}
\lambda_{\lm}\mathcal{F}(x)=\begin{cases}
1; & x\in\lm\mathcal{F}\\
0; & x\notin\lm\mathcal{F}
\end{cases},\label{eq:embedConvtoCap}
\end{equation}
and this embedding $e:\conv\to\Cap$ admits a right-adjoint $c:\Cap\rightarrow\conv$
and a left-adjoint $r:\Cap\rightarrow\conv$ defined on objects by
\begin{eqnarray}
x & \in & \lim\nolimits_{c(\lambda)}\mathcal{F\Longleftrightarrow\lambda F}(x)=1\label{eq:cCaptoConv}\\
x & \in & \lim\nolimits_{r(\lambda)}\mathcal{F\Longleftrightarrow\lambda F}(x)>0\label{eq:rCaptoConv}
\end{eqnarray}
which extend to a concrete coreflector and a concrete reflector respectively.
The functors $c$ and $r$ restrict to coreflectors and reflectors
from $\psap$ to $\pstop$, $\mathbf{\prap}$ to $\prtop$ and $\ap$
to $\Top$.

\subsection{Towards a pointfree abstraction\label{subsec:Towards-a-pointfree}}

In order to abstract the structure away to a pointfree setting, we
want to abstract $\mathbb{P}X$ to a lattice $L$ and $[0,1]$ to
a quantale $V$ as in Section \ref{subsec:V}. For now, let us assume
$V=[0,1]$ with its usual order. We first need to reinterpret $\lambda:\mathbb{FP}X\to V^{X}$
via the corresponding monotone map $\Lambda:\mathbb{FP}X\to(\mathbb{P}X)^{V}$
(giving ``convergence by levels'') defined by
\begin{equation}
\Lambda(\F)(v)=\left\{ x\in X:\lambda(\F)(x)\geq v\right\} =\lambda\F^{-1}[v,1].\label{eq:CAPLambda}
\end{equation}

Identifying $\Lambda(\F)(v)$ with $\lim_{p_{v}}\F$ for a convergence
$p_{v}$, the map $\Lambda$ corresponds to the presentation of a
$\Cap$-structure in terms of convergence tower \cite{towers}, that
is, a $\Cap$-structure on $X$ can equivalently be given by a so-called
\emph{convergence tower $\bar{p}=\{p_{\epsilon}:\epsilon\in[0,\infty]\}$
}on $X$ where each $p_{\epsilon}$ is a convergence on $X$, $p_{\infty}$
is indiscrete,
\begin{equation}
\epsilon\leq\gamma\then p_{\gamma}\leq p_{\epsilon}\label{eq:towermonotone}
\end{equation}
 and 
\begin{equation}
p_{\epsilon}=\bigvee_{\gamma>\epsilon}p_{\gamma}\label{eq:LT3}
\end{equation}
for each $\epsilon\in[0,\infty]$. The morphism condition on a convergence
tower is that $f:(X,\bar{p})\to(Y,\bar{q})$ be continuous component-wise,
that is, $f:(X,p_{\epsilon})\to(Y,q_{\epsilon})$ is continuous for
each $\epsilon$.

In our setting where $V$ is an abstraction of $[0,\infty]^{op}$
(e.g., $[0,1]$), that $p_{\infty}$ is antidiscrete becomes that
$p_{0}$ is antidiscrete, which is built in (\ref{eq:CAPLambda}).
On the other hand, (\ref{eq:towermonotone}) becomes 
\[
v\leq t\then p_{v}\leq p_{t}
\]
and corresponds to the fact that, in view of (\ref{eq:CAPLambda}),
$\Lambda(\F)$ is a decreasing function of $v$. Similarly, (\ref{eq:LT3})
becomes
\[
p_{\epsilon}=\bigvee_{\gamma<\epsilon}p_{\gamma},
\]
which in terms of $\Lambda$ is 
\[
\Lambda(\F)(\bigvee_{v\in A}v)=\bigwedge_{v\in A}\Lambda(\F)(v)
\]
for every $A\subset V$. Note that this condition also implies that
$\Lambda(\F)$ be a decreasing function of $v$.

The morphism condition becomes 
\[
\Lambda_{X}(\F)(v)\leq f^{-1}(\Lambda_{Y}(f[\F])(v)),
\]
where $\leq$ is here the order of $\mathbb{P}X$, hence inclusion.
Note that in terms of $\varphi=f^{-1}:\mathbb{P}Y\to\mathbb{P}X$,
this can be written
\[
\Lambda_{X}(\F)(v)\leq\varphi(\Lambda_{Y}(\varphi^{-1}(\F))(v)).
\]

\begin{rem}
\label{rem:classicalmorphism} Note that $\varphi=f^{-1}:\mathbb{P}Y\to\mathbb{P}X$
preserves arbitrary infima and arbitrary suprema, in particular sends
$\bot_{\mathbb{P}Y}$ to $\bot_{\mathbb{P}X}$ and $\top_{\mathbb{P}Y}$
to $\top_{\mathbb{P}X}$. It also preserves complements.
\end{rem}

\subsection{Pointfree version}

We now abstract away powersets to general lattices and $V$ is as
in Section \ref{subsec:V}. If $\mathsf{C}$ is a category of lattices,
then a $\mathsf{C}^{\Cap}$-object is a $\mathsf{C}$-object $L$
together with a monotone map 
\[
\Lambda:\mathbb{F}L\to L^{V},
\]
such that 
\begin{equation}
\Lambda(\F)(\bigvee_{v\in A}v)=\bigwedge_{v\in A}\Lambda(\F)(v)\label{eq:ordeLambda}
\end{equation}
for every $A\subset V$. In particular, taking $A=\{v_{1},v_{2}\}$
yields that for each $\F\in\mathbb{F}L$, the map $\Lambda(\F)$ is
a decreasing function of $v$. Taking $A=V$, together with (\ref{eq:density}),
yields that $\Lambda(\F)(1)=\bigwedge_{v\in V}\Lambda(\F)(v)$ and
taking $A=\emptyset$ yields that $\Lambda(\F)(0)=\top_{L}$.

A $\mathsf{C}$-morphism $\varphi:L\to L'$ between two $\mathsf{C}^{\Cap}$-objects
is a $\mathsf{C}^{\Cap}$-morphism if
\begin{equation}
\Lambda_{L'}(\F)(v)\leq\varphi(\Lambda_{L}(\varphi^{-1}(\F))(v)),\label{eq:Ccapmorphism}
\end{equation}
for every $v\in V$ (\footnote{Note that $\varphi:L\to L'$ is an abstraction of $f^{-1}:\mathbb{P}Y\to\mathbb{P}X$.}). 

Following \cite{FredetJean}, we call a category $\mathsf{C}$ of
lattices \emph{admissible }if $\mathbb{P}X$ is $\mathsf{C}$-object
for every set $X$ and there are two classes of index sets $\mathcal{I}$
and $\mathcal{J}$ such that for every pair $L,L'$ of $\mathsf{C}$-objects,
the $\mathsf{C}$-morphisms are exactly the monotone maps that preserve
all $I$-indexed infima that exist in $L$, $I\in\mathcal{I}$, and
all $J$-indexed suprema that exist in $L$, $J\in\mathcal{J}$. If
$\mathsf{C}$-morphisms additionally preserve all infima, we say that
$\mathsf{C}$ is \emph{fully admissible}.

In view of Section \ref{subsec:Towards-a-pointfree}, we have:
\begin{prop}
Let $\mathsf{C}$ be an admissible category of lattices. The association
$\mathbb{P}:\Cap\to\mathsf{C}^{\Cap}$ defined on objects by 
\[
\mathbb{P}((X,\lambda))=(\mathbb{P}X,\Lambda_{\mathbb{P}X})\text{ where }\Lambda_{\mathbb{P}X}\mathcal{F}(\alpha)=\lambda\mathcal{F}^{-1}[\alpha,1],
\]
and on morphisms by 
\[
\mathbb{P}(f)=f^{-1}:(\mathbb{P}Y,\Lambda_{\mathbb{P}Y})\to(\mathbb{P}X,\Lambda_{\mathbb{P}X})
\]
whenever $f:(X,\lambda_{X})\to(Y,\lambda_{Y})$ is a $\Cap$-morphism,
defines a contravariant functor, that is, $\mathbb{P}:\Cap\to(\mathsf{C}^{\Cap})^{\mathsf{op}}$
is a functor.
\end{prop}

\begin{rem}
\label{rem:centeredonPX} Note that if $(\mathbb{P}X,\Lambda)$ is
in the image of the functor $\mathbb{P}$, that is, if $\Lambda=\Lambda_{\mathbb{P}X}$
for a convergence approach space $(X,\lambda)$ then $x\in\Lambda(\uparrow\{x\})(v)=(\lambda\{x\}^{\uparrow})^{-1}[v,1]$
for every $v\in V$ because $\lambda\{x\}^{\uparrow}(x)=1$, that
is, 
\begin{equation}
\forall x\in X,x\in\Lambda(\uparrow\{x\})(1).\label{eq:centeredPX}
\end{equation}
\end{rem}

\section{Space of points}

The construction of a functor in the opposite direction mimics that
performed in the context of $\mathsf{C}^{\conv}$ in \cite{FredetJean}.
Several proofs are simple variations of those found in \cite{FredetJean}
and are thus omitted. To start with, we follow \cite{FredetJean}
in appealing to the categorical definition of points, that is, morphisms
whose domain is a terminal object. Singletons are terminal objects
in topological constructs. In particular, it can be verified that
$(\lbrace\ast\rbrace,\lambda)$ is a $\Cap$-terminal object. Moreover,
the corresponding $\mathsf{C}^{\Cap}$-object is initial in $\mathsf{C}^{\Cap}$
(this is an analog of \cite[Lemma 2.7]{FredetJean}, with a similar
proof):
\begin{lem}
Let $\mathsf{C}$ be an admissible category of lattices with least
and greatest elements. Then $2:=(\mathbb{P}\lbrace\ast\rbrace,\Lambda_{\mathbb{P}\lbrace\ast\rbrace})$
is an initial object in $\mathsf{C}^{\Cap}$, hence a terminal object
in $(\mathsf{C}^{\Cap})^{\mathsf{op}}$.
\end{lem}

\global\long\def\op{(\mathsf{C}{}^{\Cap})^{\mathsf{op}}}%
\global\long\def\pt{\operatorname{pt}}%

As a result, a \emph{point }of $L$ in $(\mathsf{C}^{\Cap})^{\mathsf{op}}$
is a $(\mathsf{C}^{\Cap})^{\mathsf{op}}$-morphism $\varphi:2\to L$,
that is, a $\mathsf{C}^{\Cap}$-morphism $\varphi:L\to2$. 

Following \cite[Remark 2.9]{FredetJean}, note that letting $\F=\varphi^{-1}(\{\{*\}\})\in\mathbb{F}L$,
the morphism condition (\ref{eq:Ccapmorphism}) amounts to 
\begin{equation}
\forall v\in V,\;\Lambda\F(v)\in\F,\label{eq:point}
\end{equation}
 so that points of $L$ can be identified with certain filters on
$L$ satisfying (\ref{eq:point}). What type of filters need to be
considered depends on the nature of $\mathsf{C}$-morphisms, so that
points of $L$ in $(\mathsf{C}^{\Cap})^{\mathsf{op}}$ are filters
on $L$ satisfying (\ref{eq:point}) that are additionally prime if
$\mathsf{C}$ is the category of lattices, completely prime if $\mathsf{C}$
is the category of frames, and principal filters of a join-prime element
if $\mathsf{C}$ is the category of coframes. See \cite[Remarks 2.10, 2.11, 2.13]{FredetJean}.
When $\mathsf{C=CH}$ then points are principal filters of compact
elements $\ell$ with $\ell\leq\Lambda(\uparrow\ell)(1)$.

Let $\pt_{\mathsf{C}}L$, or $\pt L$ if $\mathsf{C}$ is clear, denote
the set of all points of $L$ (in $\op$). Following a similar path
to that charted in \cite{FredetJean}, we will now endow $\pt L$
with a $\Cap$-structure. To this end, given $\ell\in L$ where $L$
is a $\mathsf{C}^{\Cap}$-object, let 
\[
\ell^{\bullet}=\{\varphi\in\pt L:\varphi(\ell)=\{*\}\}=\{\F\in\pt L:\ell\in\F\}
\]
and, given $\F\in\mathbb{FP}(\pt L)$, let 
\[
\F^{\circ}=\{\ell\in L:\ell^{\bullet}\in\F\}.
\]

The following is straightforward and can be proved along the lines
of the analogous statement \cite[Lemma 2.19]{FredetJean}:
\begin{lem}
\label{lem:bulletsandcircles} Let $\mathsf{C}$ be an admissible
category of lattices with top and bottom, and let $L$ be a $\mathsf{C}^{\Cap}$-object
and let $\ell,\ell'\in L$. 
\begin{enumerate}
\item \label{enu:bulletmonotone} if $\ell\leq\ell^{\prime}$ then $\ell^{\bullet}\subseteq{\ell^{\prime}}^{\bullet}$. 
\item \label{enu:circlemonotone} For all $\mathcal{F},\mathcal{G}\in\mathbb{FP}(\pt L)$,
if $\mathcal{F}\subset\mathcal{G}$ then $\mathcal{F}^{\circ}\subset\mathcal{G}^{\circ}$. 
\item \label{enu:bulletsup}For all $\ell_{1},\ell_{2},\cdots,\ell_{n}\in L$
$(n\in\mathbb{N})$, $\bigcup^{n}_{i=1}\ell^{\bullet}_{i}=\left(\bigvee^{n}_{i=1}\ell_{i}\right)^{\bullet}$
and if points are completely prime, then 
\[
\bigcup_{\ell\in A}\ell^{\bullet}=(\bigvee A)^{\bullet}
\]
for every $A\subset L$. 
\item \label{enu:bulletinf}For all $\ell_{1},\ell_{2},\cdots,\ell_{n}\in L$
$(n\in\mathbb{N})$, $\bigcap^{n}_{i=1}\ell^{\bullet}_{i}=\left(\bigwedge^{n}_{i=1}\ell_{i}\right)^{\bullet}$
and if points are principal (e.g., $\mathsf{C}$ is fully admissible)
then
\[
\bigcap_{\ell\in A}\ell^{\bullet}=(\bigwedge A)^{\bullet}
\]
 for every $A\subset L$. 
\item For every $\mathcal{F}\in\mathbb{FP}(\pt L)$, $\mathcal{F}^{\circ}\in\mathbb{F}L$. 
\item For every $\varphi\in\pt L$, $(\{\varphi\}^{\uparrow})^{\circ}=\varphi^{-1}\left(\left\lbrace \{*\}\right\rbrace \right)$
is the filter on $L$ representing $\varphi$. 
\item \label{enu:ptaxiominptL} $x\in\left(\Lambda_{L}(\{x\}^{\uparrow})^{\circ}(1)\right)^{\bullet}$
for every $x\in\pt L$.
\end{enumerate}
\end{lem}

Similarly to Remark \ref{rem:ptPXinCconv}, we have:
\begin{rem}
\label{rem:ptPxisX} When $\mathsf{C}$ is an admissible category
of frames or of coframes, the set $\pt(\mathbb{P}X)$ can be identified
with $X$. Indeed, under this condition, elements of $\pt(\mathbb{P}X)$
are either completely prime filters or principal filters of join-prime
elements in $\mathbb{P}X$. Either way, they are principal ultrafilters.
Remains to see that a principal ultrafilter $\{x\}^{\uparrow}\in\mathbb{FP}X$
is always an element of $\pt\mathbb{P}X$. Indeed, for every $v\in V$,
$x\in\Lambda_{\lambda}\{x\}^{\uparrow}(v)=\{t\in X:\lambda\{x\}^{\uparrow}(t)\geq v\}$
because $\lambda\{x\}^{\uparrow}(x)=1$.
\end{rem}

\begin{rem}
\label{rem:bullettopologiesonpts} Note that $\top^{\bullet}_{L}=\pt L$
and $\bot^{\bullet}_{L}=\emptyset$ so that, in view of Lemma \ref{lem:bulletsandcircles}
(\ref{enu:bulletsup}) and (\ref{enu:bulletinf}), the collection
$\{\ell^{\bullet}:\ell\in L\}$ is the collection of closed sets for
a topology on $\pt L$ whenever $\mathsf{C}$ is fully admissible,
and the collection of open sets for a topology on $\pt L$ whenever
points are completely prime, in particular when $\mathsf{C}$-morphisms
preserve arbitrary joins.
\end{rem}

\begin{defn}
Given an admissible category $\mathsf{C}$ of lattices and a $\mathsf{C}^{\Cap}$-object
$L$, $\pt L$ denotes the set of points of $L$ (in $\op$) endowed
with the $\Cap$-structure given by 
\begin{equation}
\lambda_{\pt L}\mathcal{F}(x)=\bigvee\lbrace v\in V:x\in(\Lambda\mathcal{F}^{\circ}(v))^{\bullet}\rbrace.\label{eq:lambdaptL}
\end{equation}
\end{defn}

That $\lambda_{\pt L}$ is indeed a $\Cap$-structure is easily seen:
it satisfies (\ref{eq:CAPcentered}) by Lemma \ref{lem:bulletsandcircles}(\ref{enu:ptaxiominptL})
and $\lambda_{\pt L}:\mathbb{FP}(\pt L)\to V^{X}$ is monotone because,
in view of Lemma \ref{lem:bulletsandcircles}(\ref{enu:circlemonotone})
and Lemma \ref{lem:bulletsandcircles}(\ref{enu:bulletmonotone})
and the fact that $\Lambda:\mathbb{F}L\to L^{V}$ is monotone, $(\Lambda\mathcal{F}^{\circ}(v))^{\bullet}\subset(\Lambda\mathcal{\G}^{\circ}(v))^{\bullet}$
whenever $\mathcal{F},\mathcal{G}\in\mathbb{FP}(\pt L)$ with $\F\subset\G$.

To prove that $\pt:\mathsf{C}^{\Cap}\rightarrow\Cap$ defines an adjoint
functor, we proceed as in \cite{FredetJean} showing that $\epsilon_{L}=\left(-\right)^{\bullet}:L\rightarrow\mathbb{P}\left(\pt L\right)$
are universal arrows.
\begin{prop}
If $\mathsf{C}$ is a fully admissible category of lattices, then
for every $\mathsf{C}^{\Cap}$-object $L$, the map $\epsilon_{L}:L\rightarrow\mathbb{P}(\pt L)$
defined by $\epsilon_{L}(\ell)=\ell^{\bullet}$ is a final $\mathsf{C}^{\Cap}$-morphism.
\end{prop}

\begin{proof}
Proving that $\epsilon_{L}$ is a $\mathsf{C}$-morphism can be done
exactly as in \cite[Lemma 2.23]{FredetJean}, so, we focus on proving
that $\epsilon_{L}$ is continuous, that is,
\begin{equation}
\Lambda_{\mathbb{P}\pt L}\mathcal{F}(\alpha)\leq\epsilon_{L}\Lambda\epsilon^{-1}_{L}\mathcal{F}(\alpha),\label{eqn_adjn}
\end{equation}

for every $\alpha\in V$ and $\F\in\mathbb{FP}\pt L$, that is,
\[
\{x\in\pt L:\bigvee\lbrace v\in V:x\in(\Lambda_{L}\mathcal{F}^{\circ}(v))^{\bullet}\rbrace\ge\alpha\}\subset(\Lambda_{L}\F^{\circ}(\alpha))^{\bullet}.
\]
 Let $x\in\pt L$ with $\bigvee\lbrace v\in V:x\in(\Lambda_{L}\mathcal{F}^{\circ}(v))^{\bullet}\rbrace\ge\alpha$.
In view of (\ref{eq:ordeLambda}), 
\[
\Lambda_{L}\mathcal{F}^{\circ}(\bigvee_{x\in(\Lambda_{L}\mathcal{F}^{\circ}(v))^{\bullet}}v)=\bigwedge_{x\in(\Lambda_{L}\mathcal{F}^{\circ}(v))^{\bullet}}\Lambda_{L}\F^{\circ}(v)\leq\Lambda_{L}\F^{\circ}(\alpha),
\]
and $x\in(\Lambda_{L}\mathcal{F}^{\circ}(v))^{\bullet}$ means $\Lambda_{L}\mathcal{F}^{\circ}(v)\in x$.
As $\mathsf{C}$ is fully admissible, $x$ is closed under infima
so that $\bigwedge_{x\in(\Lambda_{L}\mathcal{F}^{\circ}(v))^{\bullet}}\Lambda_{L}\F^{\circ}(v)$,
hence $\Lambda_{L}\F^{\circ}(\alpha)$, belong to $x$, that is, $x\in(\Lambda_{L}\F^{\circ}(\alpha))^{\bullet}$,
proving (\ref{eqn_adjn}). 

Conversely, if $x\in(\Lambda_{L}\F^{\circ}(\alpha))^{\bullet}$, then
$\bigvee\lbrace v\in V:x\in(\Lambda_{L}\mathcal{F}^{\circ}(v))^{\bullet}\rbrace\ge\alpha$
and thus (\ref{eqn_adjn}) is an equality and $\epsilon_{L}$ is final. 
\end{proof}

\begin{cor}
\label{cor:adjunction}If $\mathsf{C}$ is a fully admissible category
of lattices then $\pt:\mathsf{(C^{Cap})^{op}\to\Cap}$ is a right-adjoint
to $\mathbb{P}:\mathsf{Cap}\to\mathsf{(C^{Cap})^{op}}$ and thus $\Cap$
can be seen as a coreflective subcategory of $\mathsf{(C^{Cap})^{op}}$.
\end{cor}

We call a $\mathsf{C^{Cap}}$-object $(L,\Lambda)$ \emph{spatial
}if $\epsilon_{L}:L\to\mathbb{P}(\pt L)$ is an isomorphism of $\mathsf{C^{Cap}}$
and \emph{weakly spatial }if $\epsilon_{L}$ has a trivial kernel,
that is,
\begin{equation}
\ell^{\bullet}=\emptyset\iff\ell=\bot_{L}.\label{eq:WS}
\end{equation}

Assume $\pt L\neq\emptyset$. Note that in view of Lemma \ref{lem:bulletsandcircles}
(\ref{enu:bulletinf}), if $\F\in\mathbb{FP}(\pt L)$ then
\[
\F^{\natural_{\bullet}}:=(\F\cap\{\ell^{\bullet}:\ell\in L\})^{\uparrow_{\pt L}}
\]
is also a (set-theoretic) filter on $\pt L$ and $\F^{\natural_{\bullet}}\leq\F$,
but 
\begin{equation}
\lambda_{\pt L}\F=\lambda_{\pt L}\F^{\natural_{\bullet}},\label{eq:regbulletonpts}
\end{equation}
because $\F^{\circ}=(\F^{\natural_{\bullet}})^{\circ}$. Note also
that if $\G\in\mathbb{F}L$ then
\[
\G_{\bullet}:=\{g^{\bullet}:g\in\G\}^{\uparrow}
\]
 is a (possibly degenerate, but not in the weakly spatial case) set-theoretic
filter on $\pt L$ and $\G\subset(\G_{\bullet})^{\circ}$. On the
other hand, if $\F\in\mathbb{FP}(\pt L)$, $(\F^{\circ})_{\bullet}=\F^{\natural_{\bullet}}\subset\F$. 
\begin{lem}
\label{lem:grillstolift} If $(L,\Lambda)$ is a $\mathsf{C^{\Cap}}$-object
and $\F,\H\in\mathbb{FP}(\pt L)$ then 
\[
\F\#\H\then\F^{\circ}\#\H^{\circ}.
\]
Moreover, if $L$ is weakly spatial then
\[
\F^{\circ}\#\H^{\circ}\then\F^{\natural_{\bullet}}\#\H^{\natural_{\bullet}}.
\]
\end{lem}

\begin{proof}
If $\F\#\H$, $f\in\F^{\circ}$ and $h\in\H^{\circ}$ then $f^{\bullet}\in\F$
and $h^{\bullet}\in\H$ and thus $f^{\bullet}\cap h^{\bullet}=(f\wedge h)^{\bullet}\neq\emptyset$
by Lemma \ref{lem:bulletsandcircles}, so that $f\wedge h\neq0_{L}$.

Conversely, if $L$ is weakly spatial and $\F^{\circ}\#\H^{\circ}$
then $f\wedge h\neq0_{L}$ whenever $f^{\bullet}\in\F$ and $h^{\bullet}\in\H$.
By weak spatiality, $(f\wedge h)^{\bullet}=f^{\bullet}\wedge h^{\bullet}\neq\emptyset$
and thus $\F^{\natural_{\bullet}}\#\H^{\natural_{\bullet}}$.
\end{proof}

\begin{lem}
\label{lem:weaklyspatialgrill} Let $L$ be weakly spatial, $\F\in\mathbb{FP}(\pt L)$
and $\G\in\mathbb{F}L$. Then
\[
\G\#\F^{\circ}\iff\G_{\bullet}\#\F^{\natural_{\bullet}}\iff\exists\U\in\beta(\F^{\natural_{\bullet}}):\U^{\circ}\supset\G\cup\F^{\circ}.
\]
\end{lem}

\begin{proof}
If $\G\#\F^{\circ}$ then $\G_{\bullet}\#\F^{\natural_{\bullet}}$
because if $g\in\G$ and $f\in\F^{\circ}$ then $g\wedge f\neq\bot_{L}$
and 
\[
(g\wedge f)^{\bullet}=g^{\bullet}\cap f^{\bullet}\neq\emptyset
\]
 because $L$ is weakly spatial. If $\G_{\bullet}\#\F^{\natural_{\bullet}}$
there is $\U\in\beta(\G_{\bullet}\vee\F^{\natural_{\bullet}})$ and
\[
\U^{\circ}\supset(\G_{\bullet})^{\circ}\cup(\F^{\natural_{\bullet}})^{\circ}\supset\G\cup\F^{\circ}.
\]
Finally, if $\U^{\circ}\supset\G\cup\F^{\circ}$ then $\G\#\F^{\circ}$
because $\U^{\circ}\in\mathbb{F}L$. 
\end{proof}

\begin{thm}
If $\mathsf{C}$ is a fully admissible category of lattices, then
$\mathsf{C}^{\Cap}$ is a topological category over $\mathsf{C}$.
Indeed, the unique final lift $(L,\Lambda_{L})$ of a sink $(\varphi_{i}:(L_{i},\Lambda_{i})\to L)_{i\in I}$
is given by 
\[
\Lambda_{L}\F(v)=\bigwedge_{i\in I}\varphi_{i}(\Lambda_{i}(\varphi^{-1}_{i}(\F))(v)),
\]
with the convention that $\Lambda_{i}(\varphi^{-1}_{i}(\F))(v)=\top_{L_{i}}$
if $\varphi^{-1}_{i}(\F)$ is degenerate.
\end{thm}

\begin{proof}
Let $U:\mathsf{C^{Cap}}\to\mathsf{C}$ be the forgetful functor. It
is enough to show that every $U$-structure sink $(\varphi_{i}:(L_{i},\Lambda_{i})\to L)_{i\in I}$
admits a unique final lift $(L,\Lambda_{L})$. Whenever $\F\in\mathbb{F}L$,
$\varphi^{-1}_{i}(\F)$ is a filter on $L_{i}$ that is proper if
and only if $\varphi_{i}(\bot_{L_{i}})\notin\F$. Define $\Lambda_{L}:\mathbb{F}L\to L^{V}$
by
\[
\Lambda_{L}\F(v)=\bigwedge_{i\in I}\varphi_{i}(\Lambda_{i}(\varphi^{-1}_{i}(\F))(v))=\bigwedge_{i\in I,\varphi^{-1}_{i}(\F)\in\mathbb{F}L_{i}}\varphi_{i}(\Lambda_{i}(\varphi^{-1}_{i}(\F))(v)),
\]
with the convention that for improper filters, we have $\Lambda_{i}(L_{i})(v)=\top_{L_{i}}$.
It is easily verified that $\Lambda_{L}$ is monotone and, for every
$A\subset V$,
\begin{eqnarray*}
\Lambda_{L}\F(\bigvee A) & = & \bigwedge_{i\in I}\varphi_{i}(\Lambda_{i}(\varphi^{-1}_{i}(\F))(\bigvee A))\\
 & = & \bigwedge_{i\in I}\varphi_{i}(\bigwedge_{a\in A}\Lambda_{i}(\varphi^{-1}_{i}(\F))(a))\\
 & = & \bigwedge_{i\in I}\bigwedge_{a\in A}\varphi_{i}(\Lambda_{i}(\varphi^{-1}_{i}(\F))(a))
\end{eqnarray*}
because $\mathsf{C}$ is fully admissible, so that $\Lambda_{L}\F(\bigvee A)=\bigwedge_{a\in A}\Lambda_{L}\F(a)$.
Hence $(L,\Lambda_{L})$ is a $\mathsf{C^{Cap}}$-object. Moreover,
each $\varphi_{j}:(L_{j},\Lambda_{j})\to(L,\Lambda_{L})$ is a $\mathsf{C^{Cap}}$-morphism
because 
\[
\Lambda_{L}\F(v)\leq\bigwedge_{i\in I}\varphi_{i}(\Lambda_{i}(\varphi^{-1}_{i}(\F))(v))\leq\varphi_{j}(\Lambda_{j}(\varphi^{-1}_{j}(\F))(v)).
\]

To see the final lift property, let $\psi\in\mathsf{C}(L,M)$ where
$(M,\Lambda_{M})$ is a $\mathsf{C^{Cap}}$-object and suppose that
each $\psi\circ\varphi_{i}:(L_{i},\Lambda_{i})\to(M,\Lambda_{M})$
is a $\mathsf{C^{Cap}}$-morphism, that is, 
\[
\Lambda_{M}(\G)(v)\leq(\psi\circ\varphi_{i})(\Lambda_{i}((\psi\circ\varphi_{i})^{-1}(\G)(v))=\psi(\varphi_{i}(\Lambda_{i}(\varphi^{-1}_{i}(\psi^{-1}(\G))(v)))
\]
for every $\G\in\mathbb{F}M$. Thus 
\begin{eqnarray*}
\Lambda_{M}(\G)(v) & \leq & \bigwedge_{i\in I}\psi(\varphi_{i}(\Lambda_{i}(\varphi^{-1}_{i}(\psi^{-1}(\G))(v)))\\
 & \leq & \psi(\bigwedge_{i\in I}\varphi_{i}(\Lambda_{i}(\varphi^{-1}_{i}(\psi^{-1}(\G))(v)))\\
 & \leq & \psi(\Lambda_{L}(\psi^{-1}(\G))(v)),
\end{eqnarray*}
so that $\psi\in\mathsf{C^{Cap}}(L,M)$.
\end{proof}

\section{$\mathbf{\mathsf{C}^{\protect\Cap}}$ and $\mathbf{\mathsf{C}^{\protect\conv}}$}

Recall that the embedding (\ref{eq:embedConvtoCap}) of $\mathsf{Conv}$
into $\mathsf{Cap}$ is both reflective and coreflective, with reflector
$r$ and $c$ given by (\ref{eq:rCaptoConv}) and (\ref{eq:cCaptoConv})
respectively. Analogously, $(\mathsf{C^{Conv})^{op}}$ embeds in $(\mathsf{C^{Cap})^{op}}$
both reflectively and coreflectively. 

Namely, $E:(\mathsf{C^{Conv})^{op}}\to(\mathsf{C^{Cap})^{op}}$ sends
a $\mathsf{C^{Conv}}$-object $(L,\lim)$ to the $\mathsf{C^{Cap}}$-object
$E(L)=(L,\Lambda_{\lim})$ where 
\begin{equation}
\Lambda_{\lim}\F(v)=\begin{cases}
\lim\F & v>0\\
\top_{L} & v=0
\end{cases},\label{eq:embedCConvintoCCap}
\end{equation}
and a morphism $\varphi\in(\mathsf{C^{Conv})^{op}}(L',L)=(\mathsf{C^{Conv})}(L,L')$
to $E(\varphi)=\varphi$ which is a morphism of $(\mathsf{C^{Cap})}(L,L')$
because 
\[
\Lambda_{\lim_{L'}}\F(v)=\lm_{L'}\F\leq\varphi(\lm_{L}\varphi^{-1}(\F))=\varphi(\Lambda_{\lim_{L}}\varphi^{-1}(\F)(v))
\]
for every $v>0$.

This embedding has a right-adjoint $C:\mathsf{(C^{Cap})^{op}\to(C^{Conv})^{op}}$
defined on objects by $C(L,\Lambda)=(L,\lim_{C\Lambda})$ where 
\begin{equation}
\lm_{C\Lambda}\F=\Lambda\F(1)\label{eq:coreflCcaptoCconv}
\end{equation}
and a left-adjoint $R:\mathsf{(C^{Cap})^{op}\to(C^{Conv})^{op}}$
defined on objects by $R(L,\Lambda)=(L,\lim_{R\Lambda})$ where 
\begin{equation}
\lm_{R\Lambda}\F=\bigvee_{v>0}\Lambda\F(v),\label{eq:reflectCcaptoCconv}
\end{equation}
and both $R$ and $C$ are concrete functors. Moreover, the following
diagram commutes (we leave the easy verifications to the reader):

\[
\xymatrix{\conv\ar@{^{(}->}[r]_{\mathbb{P}}\ar@{_{(}->}[d]_{e} & (\Cconv)^{\mathsf{op}}\ar@{_{(}->}[d]_{E}\\
\Cap\ar@{^{(}->}[r]_{\mathbb{P}}\ar@<9pt>[u]^{r}\ar@<-7pt>[u]_{c} & \mathsf{(C^{Cap})^{op}}\ar@<9pt>[u]^{R}\ar@<-7pt>[u]_{C}
}
\]

On the other hand, it is important to note that the diagram no longer
commutes when the point-functors are added.

\[
\xymatrix{\conv\ar@{^{(}->}[r]_{\mathbb{P}}\ar@{_{(}->}[d]_{e} & (\Cconv)^{\mathsf{op}}\ar@<-7pt>[l]_{\pt}\ar@{_{(}->}[d]_{E}\\
\Cap\ar@{^{(}->}[r]_{\mathbb{P}}\ar@<9pt>[u]^{r}\ar@<-7pt>[u]_{c} & \mathsf{(C^{Cap})^{op}\ar@<-7pt>[l]_{\pt}}\ar@<9pt>[u]^{R}\ar@<-7pt>[u]_{C}
}
\]
For instance $\pt\circ R\neq r\circ\pt$: starting with a $\mathsf{C^{Cap}}$-object
$(L,\Lambda)$, $\pt\circ R(L,\Lambda)=\pt(L,\lm_{R\Lambda})$ has
$\{\varphi\in\mathsf{C}(L,2):\lim_{R\Lambda}\varphi^{-1}(\{*\})=\bigvee_{v>0}\Lambda\varphi^{-1}(\{*\})(v)\in\varphi^{-1}(\{*\})\}$
as underlying set, while $r\circ\pt(L,\Lambda)$ has underlying set
$\{\varphi\in\mathsf{C}(L,2):\forall v\in V,\Lambda\varphi^{-1}(\{*\})(v)\in\varphi^{-1}(\{*\})\}$
which is in general a strictly smaller set. 

However,
\begin{prop}
Let $\mathsf{C}$ be an admissible category of frames or of coframes.
Then
\[
\pt\circ R\circ\mathbb{P}=r,
\]
\[
\pt\circ C\circ\mathbb{P}=c,
\]
and
\[
c\circ\pt\circ E\circ\mathbb{P}=r\circ\pt\circ E\circ\mathbb{P}=\mathsf{I_{Conv}.}
\]
\end{prop}

\begin{proof}
In view of Remark \ref{rem:ptPXinCconv}, starting with a $\Cap$-object
$(X,\lambda)$, the underlying set of $\pt\circ R\circ\mathbb{P}(X,\lambda)$
and of $\pt\circ C\circ\mathbb{P}(X,\lambda)$ is $X$. Moreover,
\[
\lm_{\pt\circ R\circ\mathbb{P}X}\F=\lm_{R\Lambda_{\lambda}}\F=\{x\in X:\lambda\F(x)>0\}=\lm_{r\lambda}\F.
\]
Similarly, 
\[
\lm_{\pt\circ C\circ\mathbb{P}X}\F=\lm_{C\Lambda_{\lambda}}\F=\{x\in X:\lambda\F(x)=1\}=\lm_{c\lambda}\F.
\]

In view of Remark \ref{rem:ptPxisX}, the underlying set of $\pt\circ E\circ\mathbb{P}$
is $X$. Moreover, starting with a $\conv$-object $(X,\xi)$
\[
\lambda_{\pt\circ E\circ\mathbb{P}}\F(x)=\bigvee\{v\in V:x\in\Lambda\F(v)=\lm_{\xi}\F\}=\begin{cases}
1 & x\in\lm_{\xi}\F\\
0 & x\notin\lm_{\xi}\F
\end{cases},
\]
so that $c\circ\pt\circ E\circ\mathbb{P}(X,\xi)=r\circ\pt\circ E\circ\mathbb{P}(X,\xi)=(X,\xi)$.
\end{proof}

\section{\label{sec:Adherence}Adherence }

Let us start by inspecting the adherence function 
\begin{equation}
\adh_{\lambda}\F(x)=\bigvee_{\H\#\F}\lambda\H(x)=\bigvee_{\U\in\mathbb{U}(\F)}\lambda\U(x)\label{eq:adhCap}
\end{equation}
of a filter $\F$ on a $\mathsf{Cap}$-object $(X,\lambda)$, as defined
in, e.g., \cite{myn.APinCAP}, which applies also to subsets of $X$
via 
\[
\adh_{\lambda}A(\cdot)=\adh_{\lambda}\{A\}^{\uparrow}(\cdot),
\]
through the lens of $\mathbb{P}(X,\lambda)\in\mathsf{Ob(C^{Cap})}$,
hence seeking a function $\adh_{\mathbb{P}X}\F:V\to\mathbb{P}X$ given
by 
\begin{eqnarray*}
\adh_{\mathbb{P}X}\F(v) & = & \{x\in X:\adh_{\lambda}\F(x)\geq v\}.
\end{eqnarray*}
Since $\bigvee_{\H\#\F}\lambda\H(x)\geq v$ if and only if for every
$\epsilon<1$ there is $\H_{\epsilon}\#\F$ with $\lambda\H_{\epsilon}(x)\geq v\otimes\epsilon$,
equivalently, $x\in\Lambda_{\mathbb{P}X}\H_{\epsilon}(v\otimes\epsilon)$,
we obtain
\[
\adh_{\mathbb{P}X}\F(v)=\bigcap_{\epsilon<1}\bigcup_{\H\#\F}\Lambda_{\mathbb{P}X}\H(v\otimes\epsilon).
\]

These considerations lead to the following:
\begin{defn}
If $L$ is a $\mathsf{C^{Cap}}$-object and $\F\in\mathbb{F}L$ then
the \emph{adherence $\adh_{L}\F:V\to L$ of $\F$ }is given by
\[
\adh_{L}\F(v)=\bigwedge_{\epsilon<1}\bigvee_{\H\in\mathbb{F}L,\H\#\F}\Lambda\H(v\otimes\epsilon).
\]
\end{defn}

Note that in particular $\adh_{L}\F(0)=\top_{L}$ because $\Lambda\H(0\otimes\epsilon)=\Lambda\H(0)=\top_{L}$
for all $\H$ and $\epsilon$.
\begin{prop}
\label{prop:adhofintersection} Let $\mathsf{C}$ denote a subcategory
of that of coframes. If $\F$ and $\G$ are filters on $L\in\mathsf{Ob(C^{Cap})}$
then
\[
\adh_{L}(\F\cap\G)(v)=\adh_{L}\F(v)\vee\adh_{L}\G(v).
\]
In particular, $\F\subset\G$ implies $\adh_{L}\G\leq\adh_{L}\F$.
\end{prop}

\begin{proof}
In view of Lemma \ref{lem:grillofunion}, $\H\#(\F\cap\G)$ if and
only if $\H\#\F$ or $\H\#\G$ so that 
\[
\adh_{L}(\F\cap\G)(v)=\bigwedge_{\epsilon<1}\left(\bigvee_{\H\in\mathbb{F}L,\H\#\F}\Lambda\H(v\otimes\epsilon)\vee\bigvee_{\H\in\mathbb{F}L,\H\#\G}\Lambda\H(v\otimes\epsilon)\right)
\]
 which is $\adh_{L}\F(v)\vee\adh_{L}\G(v)$ because of the coframe
distributivity law.
\end{proof}

\begin{rem}
\label{rem:closureoperator} If $\ell\in L$ we write $\adh_{L}\ell$
for $\adh_{L}\uparrow\ell$ defining an operator $\adh_{L}:L\to L^{V}$.
In the setting of Proposition \ref{prop:adhofintersection}, since
$\uparrow(\ell\vee m)=\uparrow\ell\cap\uparrow m$ and $\ell\leq m$
if and only if $\uparrow m\subset\uparrow\ell$, this operator satisfies
\[
\ell\leq m\then\adh_{L}\ell\leq\adh_{L}m
\]
and 
\[
\adh_{L}(\ell\vee m)=\adh_{L}\ell\vee\adh_{L}m,
\]
for every $\ell,m\in L$, by Proposition \ref{prop:adhofintersection}.
Moreover, as $\bot_{L}\notin\H^{\#}$ for every proper filter and
the empty join is $\bot_{L}$, we have
\[
\adh_{L}\bot_{L}=\bot_{L}.
\]
\end{rem}

Let us note that (\ref{eq:adhCap}) allows for alternative characterizations
of objects of the full subcategories $\psap$ and $\prap$ of $\Cap$
(See, e.g., \cite{mynardmeasureCAP}). Namely, $(X,\lambda)$ is a
pseudo-approach space if and only if 
\begin{equation}
\lambda\F=\bigwedge_{\H\#\F}\adh_{\lambda}\H,\label{eq:psapincap}
\end{equation}
for every filter $\F$ on $X$, and a pre-approach space if and only
if 
\begin{equation}
\lambda\F=\bigwedge_{H\in\F^{\#}}\adh_{\lambda}H.\label{eq:prapincap}
\end{equation}

Letting $A^{(\epsilon)}=(\adh_{\lambda}A)^{-1}[\epsilon,1]$ for every
$A\subset X$, a pre-approach space $(X,\lambda)$ is an \emph{approach
space }if 
\begin{equation}
\adh_{\lambda}A(x)\geq\adh_{\lambda}A^{(\epsilon)}(x)\otimes\epsilon,\label{eq:lambdaAP}
\end{equation}
 for every $A\subset X$, every $x\in X$ and every $\epsilon\in V$.

\subsection{Adherence of image filters}

We seek to understand the effect of maps that may or may not be morphisms
in the sense of (\ref{eq:Ccapmorphism}).

Translated into our version of $\mathsf{Cap}$, the \emph{measure
of continuity of }$f:X\to Y$ where $(X,\lambda_{X})$ and $(Y,\lambda_{Y})$
are two convergence approach spaces introduced in \cite{MR3200459}
under the name \emph{default of contraction} is given by
\begin{equation}
\mathsf{m}(f)=\sup\{v\in V:\forall\F\in\mathbb{FP}X,\lambda_{Y}(f[\F])\circ f(\cdot)\geq v\otimes\lambda_{X}\F(\cdot)\}.\label{eq:mfCAP}
\end{equation}
In other words, $\lambda_{X}\F(x)\geq w$ implies $\lambda_{Y}(f[\F])(f(x))\geq v\otimes w$
for every $v$ in the defining set of $\mathsf{m}(f)$, that is,
\[
\Lambda_{X}\F(w)\subset f^{-1}(\Lambda_{Y}(f[\F])(w\otimes v)
\]
for every such $v$. These considerations lead to:
\begin{defn}
The \emph{measure of continuity $\mathsf{M}(\varphi)$ }of a map $\varphi:L\to L'$
between two $\mathsf{C^{Cap}}$-objects is given by 
\begin{equation}
\mathsf{M}(\varphi)=\bigvee\{v\in V:\forall\F\in\mathbb{F}L':\Lambda_{L'}(\F)(\cdot)\leq\varphi(\Lambda_{L}(\varphi^{-1}(\F))(\cdot\otimes v))\}.\label{eq:mphiCCap}
\end{equation}
\end{defn}

\begin{rem}
Note that in view of (\ref{eq:quantalesup}) and (\ref{eq:ordeLambda}),
letting $F=\{v\in V:\forall\F\in\mathbb{F}L':\Lambda_{L'}(\F)(\cdot)\leq\varphi(\Lambda_{L}(\varphi^{-1}(\F))(\cdot\otimes v))\}$,
\[
\Lambda_{L}(\varphi^{-1}(\F))(\cdot\otimes\mathsf{M}(\varphi))=\Lambda_{L}(\varphi^{-1}(\F))(\cdot\otimes\bigvee F)=\bigwedge_{v\in F}\Lambda_{L}(\varphi^{-1}(\F))(\cdot\otimes v)\geq\Lambda_{L'}(\F)(\cdot)
\]
for every $\F\in\mathbb{F}L'$ so that $\mathsf{M}(\varphi)\in F$
is the largest $v\in V$ satisfying the defining condition of $F$.
\end{rem}

It is plain that the morphism condition (\ref{eq:Ccapmorphism}) is
satisfied if and only if $\mathsf{M}(\varphi)=1$.
\begin{thm}
\label{thm:preserveadh} If $\mathsf{C}$ is a category of coframes,
$\varphi\in\mathsf{C}(L,L')$, and $\F\in\mathbb{F}L'$, then
\[
\adh_{L'}\F(\cdot)\leq\varphi(\adh_{L}\varphi_{\ell}(\F)(\cdot\otimes\mathsf{M}(\varphi)))\leq\varphi(\adh_{L}\varphi^{-1}(\F)(\cdot\otimes\mathsf{M}(\varphi))).
\]
\end{thm}

\begin{proof}
For each $\H\in\mathbb{F}L'$ with $\H\#\F$, $\varphi_{\ell}(\F)\#\varphi^{-1}(\H)$
by (\ref{eq:grillcoframe}). Hence, 
\[
\Lambda_{L'}\H(w\otimes\epsilon)\leq\varphi(\Lambda_{L}\varphi^{-1}(\H))(w\otimes\mathsf{M}(\varphi)\otimes\epsilon))\leq\varphi\Big(\bigvee_{\substack{\mathcal{G}\in\mathbb{F}L\\
\varphi_{\ell}(\F)\#\mathcal{G}
}
}\Lambda_{L}\mathcal{G}(w\otimes\mathsf{M}(\varphi)\otimes\epsilon)\Big).
\]
Hence, taking the supremum over all filters $\H$ meshing with $\F$
and using (\ref{eq:quantalesup}) yields
\[
\bigvee_{\H\in\mathbb{F}L,\H\#\F}\Lambda_{L'}\H(w\otimes\epsilon)\leq\varphi\Big(\bigvee_{\substack{\mathcal{G}\in\mathbb{F}L\\
\varphi_{\ell}(\F)\#\mathcal{G}
}
}\Lambda_{L}\mathcal{G}(w\otimes\mathsf{M}(\varphi)\otimes\epsilon)\Big).
\]
Finally, taking the infimum over all $\epsilon<1$ yields
\[
\bigwedge_{\epsilon<1}\bigvee_{\H\in\mathbb{F}L,\H\#\F}\Lambda_{L'}\H(w\otimes\epsilon)\leq\varphi\Big(\bigwedge_{\epsilon<1}\bigvee_{\substack{\mathcal{G}\in\mathbb{F}L\\
\varphi_{\ell}(\F)\#\mathcal{G}
}
}\Lambda_{L}\mathcal{G}(w\otimes\mathsf{M}(\varphi)\otimes\epsilon)\Big)
\]
 because $\varphi$ preserves infima, and the conclusion follows.

The rightmost inequality follows from $\varphi^{-1}(\F)\subset\varphi_{\ell}(\F)$
because $\varphi_{\ell}(\varphi(m))\leq m$. 
\end{proof}

\begin{cor}
\label{cor:preserveadhcont} If $\mathsf{C}$ is a category of coframes,
$\varphi\in\mathsf{C^{Cap}}(L,L')$, and $\F\in\mathbb{F}L'$, then
\[
\adh_{L'}\F(\cdot)\leq\varphi(\adh_{L}\varphi_{\ell}(\F)(\cdot))\leq\varphi(\adh_{L}\varphi^{-1}(\F)(\cdot)).
\]
\end{cor}

Note that Theorem \ref{thm:preserveadh} translated in $\mathsf{Cap}$
yields a more general result than \cite[(1.8)]{myn.APinCAP}, which
applies to functions that may not be morphisms. Namely:
\begin{cor}
If $(X,\lambda_{X})$ and $(Y,\lambda_{Y})$ are two convergence approach
spaces, $f:X\to Y$ and $\G\in\mathbb{FP}Y$ then
\[
\mathsf{m}(f)\otimes\adh_{\lambda_{X}}f^{-}[\G](x)\leq\adh_{\lambda_{Y}}\G(f(x))
\]
for every $x\in X$. 
\end{cor}

\subsection{Maximal filters}

Assuming the axiom of choice, every filter on a bounded lattice $L$
is contained in a maximal filter. Let us denote by $\mathbb{U}L$
the set of maximal filters on $L$ and if $\F\in\mathbb{F}L$, let
$\beta\F=\{\U\in\mathbb{U}L:\F\subset\U\}$. If $\G,\F\in\mathbb{F}L$
and $\G\#\F$ then $\{g\wedge f:g\in\G,f\in\F\}$ is a filter-base
and there is the least filter $\F\vee\G$ containing both $\F$ and
$\G$ and $\beta(\F\vee\G)=\beta\F\cap\beta\G$. 
\begin{lem}
\label{lem:adhultra} If $(L,\Lambda)$ is a $\mathsf{C^{Cap}}$-object
and $\U\in\mathbb{U}L$ then
\[
\adh_{\Lambda}\U=\Lambda\U.
\]

If $\F\in\mathbb{F}L$ then
\[
\adh_{\Lambda}\F=\bigwedge_{\epsilon<1}\bigvee_{\U\in\beta\F}\Lambda\U(v\otimes\epsilon).
\]
\end{lem}

\begin{proof}
Because $\G\#\U$ if and only if $\G\subset\U$ by maximality, 
\[
\Lambda\U(v)\leq\adh_{\Lambda}\U(v)=\bigwedge_{\epsilon<1}\bigvee_{\G\#\U}\Lambda\G(v\otimes\epsilon)\leq\bigwedge_{\epsilon<1}\Lambda\U(v\otimes\epsilon)=\Lambda\U(\bigvee_{\epsilon<1}v\otimes\epsilon)=\Lambda\U(v).
\]

The second formula follows from the observation that $\bigvee_{\G\#\F}\Lambda\G(v\otimes\epsilon)=\bigvee_{\U\in\beta\F}\Lambda\U(v\otimes\epsilon)$
because for every $\G\#\F$ there is $\U\in\beta\F$ with $\U\supset\G$.
\end{proof}

\section{Pseudo-approach, Pre-approach and Approach Spaces within $\mathsf{C^{Cap}}$}

\global\long\def\DD{\mathfrak{D}}%

Recalling Remark \ref{rem:centeredonPX}, we set out to characterize
the images of $\psap$, $\prap$ and $\ap$ under $\mathbb{P}$.

\subsection{Depth and adherence-determinedness}

The introduction of these notions is motivated by Theorem \ref{thm:PimagePsapandPrap}.

Let us fix a category $\mathsf{C}$ of lattices. If $\DD$ is a class
of collections of filters, determined by its instantiations $\DD(L)\subset\mathbb{P}(\mathbb{F}L)$
for every $\mathsf{C}$-object $L$, we say that a $\mathsf{C^{Cap}}$-object
$(L,\Lambda)$ is $\DD$\emph{-deep }if 
\begin{equation}
\Lambda(\bigcap_{\D\in\mathbb{D}}\D)=\bigwedge_{\D\in\mathbb{D}}\Lambda\D\label{eq:Ddeep}
\end{equation}
for every $\mathbb{D}\in\DD(L)$. A set $\mathbb{D}\subset\mathbb{F}L$
is \emph{finitely covering }if for every choice $(\ell_{\D})_{\D\in\mathbb{D}}\in\prod_{\D\in\mathbb{D}}\D$,
there is a finite subset $\mathbb{D}_{0}$ of $\mathbb{D}$ with $\bigvee_{\D\in\mathbb{D}_{0}}\ell_{\D}\in\bigcap_{\D\in\mathbb{D}}\D$.
Note in particular that when $L=\mathbb{P}X$ and $\F\in\mathbb{FP}X$,
then $\beta(\F)$ is finitely covering. More generally,
\begin{lem}
\label{lem:fcbeta} If a family $\mathbb{D}$ of filters on a bounded
distributive lattice $L$ is finitely covering then 
\[
\beta(\bigcap_{\D\in\mathbb{D}}\D)=\bigcup_{\D\in\mathbb{D}}\beta\D,
\]
and conversely if $L$ is a Boolean Algebra.
\end{lem}

\begin{proof}
The inclusion $\bigcup_{\D\in\mathbb{D}}\beta\D\subset\beta(\bigcap_{\D\in\mathbb{D}}\D)$
is always true. Suppose $\mathbb{D}$ is finitely covering and let
$\F=\bigcap_{\D\in\mathbb{D}}\D$. If $\U\in\beta\F\setminus\bigcup_{\D\in\mathbb{D}}\beta\D$
then for each $\D\in\mathbb{D}$ there is $\ell_{\D}\in\D\setminus\U$.
There is a finite subset $\mathbb{D}_{0}$ of $\mathbb{D}$ with $\bigvee_{\D\in\mathbb{D}_{0}}\ell_{\D}\in\F\subset\U$.
Since $L$ is a bounded distributive lattice, maximal filters, in
particular $\U$, are prime. Hence, there is $\ell_{\D}\in\U$, which
is a contradiction. 

Conversely, assume $\beta(\bigcap_{\D\in\mathbb{D}}\D)=\bigcup_{\D\in\mathbb{D}}\beta\D$
and that $L$ is a Boolean algebra. Pick $\ell_{\D}\in\D$ for each
$\D\in\mathbb{D}$. Let $\F=\bigcap_{\D\in\mathbb{D}}\D$. If $\bigvee_{\D\in\mathbb{D}_{0}}\ell_{\D}\notin\F$
whenever $\mathbb{D}_{0}$ is a finite subset of $\mathbb{D}$, then
$(\bigvee_{\D\in\mathbb{D}_{0}}\ell_{\D})^{*}\in\F^{\#}$ and $\{(\bigvee_{\D\in\mathbb{D}_{0}}\ell_{\D})^{*}\wedge f:f\in\F,\mathbb{D}_{0}\in[\mathbb{D}]^{<\infty}\}$
is a filter-base and there is $\U\in\beta\F$ with $(\bigvee_{\D\in\mathbb{D}_{0}}\ell_{\D})^{*}=\bigwedge_{\D\in\mathbb{D}_{0}}\ell^{*}_{\D}\in\U$
for every $\mathbb{D}_{0}\in[\mathbb{D}]^{<\infty}$, so that $\ell^{*}_{\D}\in\U$
for all $\D\in\mathbb{D}$ and thus $\U\notin\bigcup_{\D\in\mathbb{D}}\beta\D=\beta\F$,
which is a contradiction.
\end{proof}

We say that $(L,\Lambda)$ is \emph{deep }if it is $\mathbb{P}(\mathbb{F}L)$-deep,
and \emph{fc-deep }if it is deep with respect to the class $\DD$
of finitely covering\emph{ }sets of filters.

On the other hand, if $\mathbb{D}$ is a class of filters determined
by its instantiations $\mathbb{D}(L)\in\mathbb{P}(\mathbb{F}L$),
we say that a $\mathsf{C^{Cap}}$-object $(L,\Lambda)$ is $\mathbb{D}$-\emph{adherence
determined }if 
\begin{equation}
\Lambda(\F)=\bigwedge_{\mathbb{D}(L)\ni\D\#\F}\adh_{\Lambda}\D,\label{eq:pfDadhdeter}
\end{equation}
for every $\F\in\mathbb{F}L$. Given $\mathsf{C}$ and a class $\mathbb{D}$
of filters and a class $\DD$ of collections of filters as above,
let us denote by $\mathsf{Adh_{\mathbb{D}}}$ and by $\mathsf{Deep_{\DD}}$
the full subcategories of $\mathsf{C^{\Cap}}$ whose objects are $\mathbb{D}$-adherence
determined and $\DD$-deep respectively.

Let $\mathbb{F}$ denote the class of all filters and $\mathbb{F}_{0}$
the class of principal filters, so that $\mathbb{F}_{0}L$ is the
set of principal filters on $L$. 
\begin{prop}
\label{prop:adhtodeep} Let $(L,\Lambda)$ be a $\mathsf{C^{Cap}}$-object.
\begin{enumerate}
\item \label{enu:frame} If $L$ is frame that is $\mathbb{F}_{0}$-adherence
determined, then $L$ is deep.
\item \label{enu:Boolean} If a bounded distributive lattice is $\mathbb{F}$-adherence
determined then it is fc-deep, and conversely in a Boolean algebra.
\end{enumerate}
\end{prop}

\begin{proof}
(\ref{enu:frame}). In a frame, $(\bigcap_{\D\in\mathbb{D}}\D)^{\#}=\bigcup_{\D\in\mathbb{D}}\D^{\#}$,
so that when $L$ is $\mathbb{F}_{0}$-determined, then
\[
\Lambda(\bigcap_{\D\in\mathbb{D}}\D)=\bigwedge_{\ell\in(\bigcap_{\D\in\mathbb{D}}\D)^{\#}}\adh_{\Lambda}\ell=\bigwedge_{\ell\in\bigcup_{\D\in\mathbb{D}}\D^{\#}}\adh_{\Lambda}\ell=\bigwedge_{\D\in\mathbb{D}}\bigwedge_{\ell\in\D^{\#}}\adh_{\Lambda}\ell=\bigwedge_{\D\in\mathbb{D}}\Lambda\D,
\]
and $L$ is deep. 

(\ref{enu:Boolean}). In view of Lemma \ref{lem:adhultra}, $(L,\Lambda)$
is $\mathbb{F}$-adherence determined if and only if
\begin{equation}
\Lambda(\F)=\bigwedge_{\H\#\F}\adh_{\Lambda}\H=\bigwedge_{\U\in\beta\F}\Lambda(\U).\label{eq:Fadhdet}
\end{equation}
 If $\mathbb{D}$ is finitely covering and $(L,\Lambda)$ is $\mathbb{F}$-adherence
determined then, letting $\F=\bigcap_{\D\in\mathbb{D}}\D$, Lemma
\ref{lem:fcbeta} applies to the effect that
\[
\Lambda\F=\bigwedge_{\U\in\beta\F=\bigcup_{\D\in\mathbb{D}}\beta\D}\Lambda\U=\bigwedge_{\D\in\mathbb{D}}\bigwedge_{\U\in\beta\D}\Lambda\U=\bigwedge_{\D\in\mathbb{D}}\Lambda\D
\]
and $(L,\Lambda)$ is fc-deep.

Conversely, if $(L,\Lambda)$ is fc-deep and $L$ is a Boolean algebra,
then each set $\beta\F$ is finitely covering and thus (\ref{eq:Fadhdet}).
\end{proof}

\begin{rem}
When $L=\mathbb{P}X$, the converse to Proposition \ref{prop:adhtodeep}
(\ref{enu:frame}) is true. Indeed, if $(\mathbb{P}X,\Lambda)$ is
deep, $x\in X$ and $w\in V$, let $\N_{w}(x)=\bigcap\{\D\in\mathbb{FP}X:x\in\Lambda\D(w)\}$.
It is easily checked that
\[
x\in\Lambda\F(w)\iff\N_{w}(x)\subset\F\iff\F^{\#}\subset\N_{w}(x)^{\#}.
\]
Hence, if 
\[
x\in\bigcap_{A\in\F^{\#}}\adh_{\Lambda}A(v)=\bigcap_{A\in\F^{\#}}\bigcap_{\epsilon<1}\bigcup_{\H\#A}\Lambda\H(v\otimes\epsilon)
\]
then for every $A\in\F^{\#}$ and $\epsilon<1$, there is $\H\#A$
with $x\in\Lambda\H(v\otimes\epsilon)$, so that $A\in\H^{\#}\subset(\N_{v\otimes\epsilon}(x))^{\#}$.
Thus $\F^{\#}\subset(\N_{v\otimes\epsilon}(x))^{\#}$ for every $\epsilon<1$,
that is, $x\in\bigcap_{\epsilon<1}\Lambda\F(v\otimes\epsilon)=\Lambda\F(v)$.

In contrast a deep $\mathsf{CH}^{\Cap}$-object need not be $\mathbb{F}_{0}$-adherence
determined, in fact not even $\mathbb{F}$-adherence determined, even
if weakly spatial. Consider for instance, the Heyting algebra $L=\{0,d,1\}$
with $0<d<1$ and $d^{*}=0$. There are only two proper filters, $\uparrow1$
and $\uparrow d$ and $\uparrow1\subset\uparrow d$. Let $\Lambda$
be given by $\Lambda(\uparrow1)(v)=d$ and $\Lambda(\uparrow d)(v)=1$
for every $v>0$, which defines a weakly spatial $\mathsf{CH}^{\Cap}$-object.
The only point is $\uparrow d$ so that $\pt L$ is a singleton and
is thus trivially an approach space. However, $L$ fails to be even
$\mathbb{F}$-adherence determined because $\Lambda(\uparrow1)(1)=d<1=\bigwedge_{\H\#\uparrow1}\adh_{\Lambda}\H(1)$.
\end{rem}

\subsection{\label{subsec:imagesunderP}Images of $\protect\Cap$, $\protect\psap$,
$\protect\prap$, and $\protect\ap$ under $\mathbb{P}$}
\begin{thm}
\label{thm:PimagePsapandPrap} A $\mathsf{C^{Cap}}$-object $(\mathbb{P}X,\Lambda)$
is in the image under $\mathbb{P}$ of
\begin{enumerate}
\item \label{enu:Pofcap}$\Cap$ if and only if (\ref{eq:centeredPX});
\item \label{enu:PofPsap} $\psap$ if and only if (\ref{eq:centeredPX})
and it is fc-deep, if and only if (\ref{eq:centeredPX}) and it is
$\mathbb{F}$-adherence determined;
\item \label{enu:PofPrap} $\prap$ if and only if (\ref{eq:centeredPX})
and it is deep, if and only if (\ref{eq:centeredPX}) and it is $\mathbb{F}_{0}$-adherence
determined.
\end{enumerate}
\end{thm}

\begin{proof}
(\ref{enu:Pofcap}). It is clear that the image of a $\Cap$-object
under $\mathbb{P}$ satisfies (\ref{eq:centeredPX}). Assume conversely
that $L=(\mathbb{P}X,\Lambda)$ satisfies (\ref{eq:centeredPX}).
Then $\mathbb{P}(\pt L,\lambda_{\pt L})=(\mathbb{P}X,\Lambda)$. Indeed,
$\pt L=X$ and $\F^{\circ}=\F$ so that
\begin{eqnarray*}
\Lambda_{\mathbb{P}\pt L}\F(v) & = & \{x\in X:\lambda_{\pt L}\F(x)\geq v\}\\
 & = & \{x\in X:\bigvee\{t\in V:\Lambda\F^{\circ}(t)\in\uparrow x\}\geq v\}.\\
 & = & \{x\in X:x\in\Lambda\F(v)\}=\Lambda\F(v).
\end{eqnarray*}

Regarding (\ref{enu:PofPsap}) and (\ref{enu:PofPrap}), note that
the equivalences between fc-deep and $\mathbb{F}$-adherence determined,
and between deep and $\mathbb{F}_{0}$-adherence determined, follow
from Proposition \ref{prop:adhtodeep}. 

(\ref{enu:PofPsap}). If $(X,\lambda)$ is a pseudo-approach space
then $(\mathbb{P}X,\Lambda)$ is $\mathbb{F}$-adherence determined.
Indeed, 
\begin{eqnarray*}
\Lambda\F(v) & = & \{x\in X:\lambda\F(x)=\bigwedge_{\H\#\F}\adh_{\lambda}\H(x)\geq v\}\\
 & = & \bigcap_{\H\#\F}\{\adh_{\lambda}\H\geq v\}=\bigwedge\adh_{\Lambda}\H(v),
\end{eqnarray*}
 because $\adh_{\Lambda}\H(v)=\{\adh_{\lambda}\H\geq v\}$ as noted
in Section \ref{sec:Adherence}. On the other hand, if $(\mathbb{P}X,\Lambda)$
is $\mathbb{F}$-adherence determined and (\ref{eq:centeredPX}),
then $(\mathbb{P}X,\Lambda)=\mathbb{P}(X,\lambda)$ where $\lambda\F(x)=\bigvee\{v\in V:x\in\Lambda\F(v)\}$.
To see that $\lambda$ is a pseudo-approach structure, first note
that it satisfies (\ref{eq:CAPcentered}) because of (\ref{eq:centeredPX}).
Moreover, $\adh_{\lambda}\H(x)=\bigvee_{\G\#\H}\bigvee_{x\in\Lambda\G(v)}v$,
so that $\adh_{\lambda}\H(t)\geq v$ if for every $\epsilon<1$, there
is $\G_{\epsilon}\#\H$ and there is $w_{\epsilon}\in V$ with $t\in\Lambda\G_{\epsilon}(w_{\epsilon})$
with $w_{\epsilon}\geq v\otimes\epsilon$, so that $t\in\Lambda\G_{\epsilon}(v\otimes\epsilon)$.
In other words, 
\begin{equation}
\{\adh_{\lambda}\H\geq v\}=\adh_{\Lambda}\H(v).\label{eq:auxadhlevel}
\end{equation}
Hence,
\begin{eqnarray*}
\lambda\F(x) & = & \bigvee\{v\in V:x\in\Lambda\F(v)=\bigwedge_{\H\#\F}\adh_{\Lambda}\H(v)\}\\
 & = & \bigvee\{v\in V:\H\#\F\then\adh_{\lambda}\H(x)\geq v\}\}\\
 & = & \bigwedge_{\H\#\F}\adh_{\lambda}\H(x).
\end{eqnarray*}

(\ref{enu:PofPrap}). Similarly, if $(X,\lambda)$ is a pre-approach
space then $(\mathbb{P}X,\Lambda)$ is deep because 
\begin{eqnarray*}
\Lambda\Big(\bigcap_{\D\in\mathbb{D}}\mathcal{\D}\Big)(v) & = & \{x\in X:\lambda\Big(\bigcap_{\D\in\mathbb{D}}\mathcal{\D}\Big)(x)=\bigwedge_{\D\in\mathbb{D}}\lambda\D(x)\geq v\}\\
 & = & \bigcap_{\D\in\mathbb{D}}\{\lambda\D\geq v\}=\bigwedge_{\D\in\mathbb{D}}\Lambda\D(v).
\end{eqnarray*}

On the other hand, if $(\mathbb{P}X,\Lambda)$ is deep and (\ref{eq:centeredPX}),
then $(\mathbb{P}X,\Lambda)=\mathbb{P}(X,\lambda)$ where $\lambda\F(x)=\bigvee\{v\in V:x\in\Lambda\F(v)\}$
satisfies (\ref{eq:CAPcentered}) because of (\ref{eq:centeredPX})
and
\begin{eqnarray*}
\lambda\Big(\bigcap_{\D\in\mathbb{D}}\mathcal{\D}\Big)(x) & = & \bigvee\{v\in V:x\in\Lambda\Big(\bigcap_{\D\in\mathbb{D}}\mathcal{\D}\Big)(v)=\bigcap_{\D\in\mathbb{D}}\Lambda\D(v)\}\\
 & = & \bigwedge_{\D\in\mathbb{D}}\lambda\D(x).
\end{eqnarray*}
\end{proof}

On the other hand, recall that an approach space is a pre-approach
space satisfying (\ref{eq:lambdaAP}).
\begin{prop}
\label{prop:PimageAp} A $\mathsf{C^{Cap}}$-object $(\mathbb{P}X,\Lambda)$
is in the image of $\ap$ under $\mathbb{P}$ if and only if it is
in the image of $\prap$ and
\begin{equation}
\adh_{\Lambda}\ell\left(v\otimes\epsilon\right)\geq\adh_{\Lambda}\Big(\left(\adh_{\Lambda}\ell\right)\left(\epsilon\right)\Big)\left(v\right)\label{eq:pfAP}
\end{equation}
for every $\ell\in\mathbb{P}X$ and $v,\epsilon\in V$.
\end{prop}

\begin{proof}
If $(X,\lambda)$ is an approach space, it is in particular a pre-approach
space, hence $(\mathbb{P}X,\Lambda_{\mathbb{P}X})$ is deep and centered.
Moreover, as $\adh_{\lambda}A(x)\geq\adh_{\lambda}A^{(\epsilon)}(x)\otimes\epsilon$
for every $A\in\mathbb{P}X$ and $A^{(\epsilon)}=\{\adh A\geq\epsilon\}=\adh_{\Lambda}A(\epsilon)$,
we obtain that if $x\in\adh_{\Lambda}\Big(\left(\adh_{\Lambda}A\right)\left(\epsilon\right)\Big)\left(v\right)$
then $\adh_{\lambda}(A^{(\epsilon)})\geq v$ thus 
\[
\adh_{\lambda}A(x)\geq\adh_{\lambda}A^{(\epsilon)}(x)\otimes\epsilon\geq v\otimes\epsilon,
\]
so that $x\in\adh_{\Lambda}A(v\otimes\epsilon)$ and (\ref{eq:pfAP})
is satisfied.

Moreover, if $(\mathbb{P}X,\Lambda)$ is centered and deep, we have
$\Lambda=\Lambda_{\mathbb{P}X}$ for the pre-approach structure $\lambda$
given by $\lambda\F(x)=\bigvee\{v\in V:x\in\Lambda\F(v)\}$, by Theorem
\ref{thm:PimagePsapandPrap}. If it satisfies (\ref{eq:pfAP}) as
well, then in view of (\ref{eq:auxadhlevel}), (\ref{eq:pfAP}) rephrases
as 
\[
\{\adh_{\lambda}A\geq v\otimes\epsilon\}\supset\{\adh_{\lambda}\Big(A^{(\epsilon)}\Big)\geq v\},
\]
yielding (\ref{eq:lambdaAP}).
\end{proof}

\subsection{Adherence in $(L,\Lambda)$ vs adherence in $\protect\pt L$}

The correspondence (\ref{eq:auxadhlevel}) between $\adh_{\lambda}$
and $\adh_{\Lambda_{\mathbb{P}X}}$ extends partially to $\adh_{\pt L}$
and $\adh_{\Lambda}$ in $\mathsf{CH^{Cap}}$:
\begin{prop}
\label{prop:adhdual-1} Let $\mathsf{C}$ be such that points of $\mathsf{C^{Cap}}$
are principal. Let $(L,\Lambda)$ be a $\mathsf{C^{Cap}}$-object.
Let $\H\in\mathbb{FP}(\pt L)$, $x\in\pt L$, and $v\in V$. Then
\[
\adh_{\pt L}\H(x)\geq v\then\adh_{\Lambda}\H^{\circ}(v)\in x.
\]
If moreover $L$ is weakly spatial and points of $\mathsf{C^{Cap}}$
are principal filters of completely prime elements (e.g., for $\mathsf{C=CH}$),
then
\[
\adh_{\Lambda}\H^{\circ}(v)\in x\then\adh_{\pt L}\H^{\natural_{\bullet}}(x)\geq v.
\]
\end{prop}

\begin{proof}
Assume that 
\[
\adh_{\pt L}\H(x)=\bigvee_{\U\in\beta(\H)}\bigvee_{\Lambda\U^{\circ}(w)\in x}w\geq v.
\]
For every $\epsilon<1$, there is $\U\in\beta(\H)$ with $\Lambda\U^{\circ}(v\otimes\epsilon)\in x$.
As $\U^{\circ}\supset\H^{\circ}$, $\bigvee_{\L\#\H^{\circ}}\Lambda\L(v\otimes\epsilon)\in x$
and $x$ is principal, so that $\adh_{\Lambda}\H^{\circ}(v)=\bigwedge_{\epsilon<1}\bigvee_{\L\#\H^{\circ}}\Lambda\L(v\otimes\epsilon)\in x$.

Suppose $\adh_{\Lambda}\H^{\circ}(v)\in x$ and $x$ is a principal
completely prime filter. Since $\adh_{\Lambda}\H^{\circ}(v)=\bigwedge_{\epsilon<1}\bigvee_{\L\#\H^{\circ}}\Lambda\L(v\otimes\epsilon)$,
$\bigvee_{\L\#\H^{\circ}}\Lambda\L(v\otimes\epsilon)\geq\bigwedge x$
for every $\epsilon<1$. Hence, by complete primeness of $\bigwedge x$,
there is $\L_{\epsilon}\#\H^{\circ}$ with $\Lambda\L_{\epsilon}(v\otimes\epsilon)\in x$.
In view of Lemma \ref{lem:weaklyspatialgrill}, there is $\U_{\epsilon}\in\beta(\H^{\natural_{\bullet}})$
with $\U^{\circ}_{\epsilon}\supset\L_{\epsilon}\cup\H^{\circ}$ and
thus $\Lambda\U^{\circ}_{\epsilon}(v\otimes\epsilon)\in x.$ Hence
\[
\adh_{\pt L}\H^{\natural_{\bullet}}(x)=\bigvee_{\U\in\beta(\H^{\natural_{\bullet}})}\lambda_{\pt L}\U(x)=\bigvee_{\U\in\beta(\H^{\natural_{\bullet}})}\bigvee_{\Lambda\U^{\circ}(w)\in x}w\geq\bigvee_{\epsilon<1}v\otimes\epsilon=v.
\]
\end{proof}

\begin{cor}
\label{cor:adhAepsilon} Let $\mathsf{C}$ be such that points of
$\mathsf{C^{Cap}}$ are principal. Let $(L,\Lambda)$ be a $\mathsf{C^{Cap}}$-object.
If $A\subset\pt L$ then $(\{A\}^{\uparrow})^{\circ}$ is principal.
Given $\epsilon\in V$, let $a_{\epsilon}:=\bigwedge(\{A^{(\epsilon)}\}^{\uparrow})^{\circ}$
and $a=:\bigwedge(\{A\}^{\uparrow})^{\circ}$. Then 
\begin{enumerate}
\item 
\begin{equation}
a_{\epsilon}\leq\adh_{\Lambda}a(\epsilon);\label{eq:aepsilon}
\end{equation}
\item If moreover, $(L,\Lambda)$ satisfies (\ref{eq:pfAP}) then 
\begin{equation}
x\in(\adh_{\Lambda}a(\adh_{\pt L}A^{(\epsilon)}(x)\otimes\epsilon))^{\bullet},\label{eq:quasidiagonal}
\end{equation}
for every $x\in\pt L$;
\item If moreover $L$ is weakly spatial and points of $\mathsf{C^{\Cap}}$
are principal completely prime filters, then 
\begin{equation}
\adh_{\pt L}a^{\bullet}(x)\geq\adh_{\pt L}A^{(\epsilon)}(x)\otimes\epsilon.\label{eq:diagonalifbulletisset}
\end{equation}
\end{enumerate}
\end{cor}

\begin{proof}
(1). Since $\ell\in(\{A\}^{\uparrow})^{\circ}$ if and only if $A\subset\ell^{\bullet}$,
we have, in view of Lemma \ref{lem:bulletsandcircles}(\ref{enu:bulletinf}),
\[
A\subset\bigwedge_{\ell\in(\{A\}^{\uparrow})^{\circ}}\ell^{\bullet}=(\bigwedge_{\ell\in(\{A\}^{\uparrow})^{\circ}}\ell)^{\bullet}
\]
so that $a:=\bigwedge_{\ell\in(\{A\}^{\uparrow})^{\circ}}\ell\in(\{A\}^{\uparrow})^{\circ}$,
and $(\{A\}^{\uparrow})^{\circ}$ is thus principal. Recall that $A^{(\epsilon)}=\{t\in\pt L:\adh_{\pt L}A(t)\geq\epsilon\}$.
For each $t\in A^{(\epsilon)}$, Proposition \ref{prop:adhdual-1}
applies to the effect that 
\[
\adh_{\Lambda}a(\epsilon)=\adh_{\Lambda}(\{A\}^{\uparrow})^{\circ}(\epsilon)\in t,
\]
that is, $t\in(\adh_{\Lambda}a(\epsilon))^{\bullet}$. In other words,
$A^{(\epsilon)}\subset(\adh_{\Lambda}a(\epsilon))^{\bullet}$, that
is, $\adh_{\Lambda}a(\epsilon)\in(\{A^{(\epsilon)}\}^{\uparrow})^{\circ}=\uparrow a_{\epsilon}$,
which yields the conclusion.

(2). Let $x\in\pt L$ and $v=\adh_{\pt L}A^{(\epsilon)}(x)$. By Proposition
\ref{prop:adhdual-1},
\[
\adh_{\Lambda}a_{\epsilon}(v)=\adh_{\Lambda}(\{A^{(\epsilon)}\}^{\uparrow})^{\circ}(v)\in x
\]
and, as $\adh_{\Lambda}:L\to L^{V}$ is monotone, its application
to (\ref{eq:aepsilon}) yields
\[
\adh_{\Lambda}a_{\epsilon}(v)\leq\adh_{\Lambda}(\adh_{\Lambda}a(\epsilon))(v)\leq\adh_{\Lambda}a(v\otimes\epsilon),
\]
where the second inequality follows from (\ref{eq:pfAP}). Since the
left-hand side is in the filter $x$, so is the right-hand side and
the conclusion follows.

(3). (\ref{eq:quasidiagonal}) rephrases as $\adh_{\Lambda}(\{A\}^{\uparrow})^{\circ}(v)\in x$
for $v=\adh_{\pt L}A^{(\epsilon)}(x)\otimes\epsilon$. Under the new
assumptions, the second part of Proposition \ref{prop:adhdual-1}
applies to $\H=\{A\}^{\uparrow}$ to the effect that
\[
\adh_{\pt L}((\{A\}^{\uparrow})^{\circ})_{\bullet}(x)=\adh_{\pt L}a^{\bullet}(x)\geq\adh_{\pt L}A^{(\epsilon)}(x)\otimes\epsilon.
\]
\end{proof}

\subsection{Depth of spaces of points}

In view of Theorem \ref{thm:PimagePsapandPrap}, $\DD$-deep $\mathsf{C^{\Cap}}$-objects
naturally generalize $\prap$ and $\psap$ objects. Moreover, we can
consider $\DD$-depth in $\Cap$ as well. Indeed, given a class $\DD$
of collections of filters, a convergence approach space $(X,\lambda)$
is called $\DD$-\emph{deep} if 
\[
\lambda(\bigcap_{\D\in\mathbb{D}}\D)=\bigwedge_{\D\in\mathbb{D}}\lambda\D,
\]
for every $\mathbb{D}\in\DD(\mathbb{P}X).$ Pre-approach spaces are
exactly $\DD$-deep convergence approach spaces where $\DD(L)=\mathbb{P}(\mathbb{F}L)$
while pseudo-approach spaces are exactly $\DD$-deep convergence approach
spaces where $\DD$ is the class of finitely covering sets of filters
(\footnote{Indeed, $\lambda(\F)=\bigwedge_{\U\in\beta\F}\lambda\U$ characterizes
pseudo-approach spaces and sets of the form $\beta\F$ are finitely
covering in $\mathbb{P}X$, so that a $\DD$-deep convergence approach
is pseudo-approach. On the other hand if $(X,\lambda)$ is pseudo-approach
and $\mathbb{D}$ is finitely covering then letting $\F=\bigcap_{\D\in\mathbb{D}}\D$,
Lemma \ref{lem:fcbeta} applies to the effect that
\[
\lambda(\F)=\bigwedge_{\U\in\beta\F=\bigcup_{\D\in\mathbb{D}}\beta\D}\lambda\U=\bigwedge_{\D\in\mathbb{D}}\bigwedge_{\U\in\beta\D}\lambda\U=\bigwedge_{\D\in\mathbb{D}}\lambda\D.
\]
}). 

Let us say that $\DD$ is $\pt$-\emph{compatible }if $\{\D^{\circ}:\D\in\mathbb{D}\}\in\DD(L)$
whenever $\mathbb{D}\in\DD(\mathbb{P}(\pt L))$. 
\begin{lem}
\label{lem:ptcomptible} The class $\DD$ given by $\DD(L)=\mathbb{P}(\mathbb{F}L)$
and the class of finitely covering sets of filters are both $\pt$-compatible.
\end{lem}

\begin{proof}
It is obvious for $\DD(L)=\mathbb{P}(\mathbb{F}L)$. Let now $\mathbb{D}\subset\mathbb{FP}(\pt L)$
be finitely covering and let $(\ell_{\D})_{\D\in\mathbb{D}}\in\prod_{\D\in\mathbb{D}}\D^{\circ}$,
that is, $\ell^{\bullet}_{\D}\in\D$ for every $\D\in\mathbb{D}$.
Because $\mathbb{D}$ is finitely covering, there is a finite subset
$\mathbb{D}_{0}$ of $\mathbb{D}$ with $\bigvee_{\D\in\mathbb{D}_{0}}(\ell_{\D})^{\bullet}\in\bigcap_{\D\in\mathbb{D}}\D$.
In view of Lemma \ref{lem:bulletsandcircles}, 
\[
\bigvee_{\D\in\mathbb{D}_{0}}(\ell_{\D})^{\bullet}=(\bigvee_{\D\in\mathbb{D}_{0}}\ell_{\D})^{\bullet}\in(\bigcap_{\D\in\mathbb{D}}\D)^{\circ}=\bigcap_{\D\in\mathbb{D}}\D^{\circ}
\]
 and $\{\D^{\circ}:\D\in\mathbb{D}\}$ is finitely covering.
\end{proof}

Remarkably, when $\DD$ is $\pt$-compatible and $\mathsf{C}$ is
fully admissible, $\DD$-depth is preserved by the point-functor:
\begin{thm}
\label{thm:deptofpoints} Let $\mathsf{C}$ be fully admissible and
let $\DD$ be a $\pt$-compatible class of sets of filters. The space
of points of a $\DD$-deep $\mathsf{C^{\Cap}}$-object is a $\DD$-deep
convergence approach space. 
\end{thm}

\begin{proof}
Let $(L,\Lambda)$ be $\DD$-deep. Let $\mathbb{D}\in\DD(\mathbb{P}\pt L)$.
As $\lambda_{\pt L}(\bigcap_{\D\in\mathbb{D}}\D)(x)=\bigvee\{v\in V:\Lambda(\bigcap_{\D\in\mathbb{D}}\D)^{\circ}(v)\in x\}$
and $(\bigcap_{\D\in\mathbb{D}}\D)^{\circ}=\bigcap_{\D\in\mathbb{D}}\D^{\circ}$,
that $\Lambda$ be $\DD$-deep yields $\Lambda(\bigcap_{\D\in\mathbb{D}}\D)^{\circ}(v)=\bigwedge_{\D\in\mathbb{D}}\Lambda\D^{\circ}(v)$.
If $\mathsf{C}$ is fully admissible, then points are principal filters
and thus $\bigwedge_{\D\in\mathbb{D}}\Lambda\D^{\circ}(v)\in x$ if
and only if $\Lambda\D^{\circ}(v)\in x$ for every $\D\in\mathbb{D}$.
In other words,
\[
\{v\in V:\Lambda(\bigcap_{\D\in\mathbb{D}}\D)^{\circ}(v)\in x\}=\bigcap_{\D\in\mathbb{D}}\{v\in V:\Lambda\D{}^{\circ}(v)\in x\}.
\]

Moreover, for each $\D\in\mathbb{D}$, in view of (\ref{eq:ordeLambda}),
\[
\Lambda\D^{\circ}(\bigvee\{v\in V:\Lambda\D{}^{\circ}(v)\in x\})=\bigwedge_{\{v\in V:\Lambda\D{}^{\circ}(v)\in x\}}\Lambda\D^{\circ}(v)
\]
belongs to $x$ because $x$ is principal, so that each set $\{v\in V:\Lambda\D{}^{\circ}(v)\in x\}$
has a greatest element. Since these sets are down-closed, the infimum
of these greatest elements is in $\bigcap_{\D\in\mathbb{D}}\{v\in V:\Lambda\D{}^{\circ}(v)\in x\}$,
so that 
\[
\bigwedge_{\D\in\mathbb{D}}\bigvee\{v\in V:\Lambda\D{}^{\circ}(v)\in x\}\leq\bigvee\bigcap_{\D\in\mathbb{D}}\{v\in V:\Lambda\D{}^{\circ}(v)\in x\}
\]
while the reverse inequality is trivial. Thus
\begin{eqnarray*}
\lambda_{\pt L}(\bigcap_{\D\in\mathbb{D}}\D)(x) & = & \bigvee\bigcap_{\D\in\mathbb{D}}\{v\in V:\Lambda\D{}^{\circ}(v)\in x\}\\
 & = & \bigwedge_{\D\in\mathbb{D}}\bigvee\{v\in V:\Lambda\D{}^{\circ}(v)\in x\}\\
 & = & \bigwedge_{\D\in\mathbb{D}}\lambda_{\pt L}\D(x),
\end{eqnarray*}
that is, $(\pt L,\lambda_{\pt L})$ is a $\DD$-deep convergence approach
space. 
\end{proof}

\subsection{Alternatives for pointfree versions of $\protect\psap$, $\protect\prap$,
and $\protect\ap$.}

Let us assume that $\mathsf{C}$ is a fully admissible category of
lattices.

In view of the results of Section \ref{subsec:imagesunderP}, we have
two natural paths to consider pointfree versions of $\psap$, $\prap$,
and $\ap$: via depth, or via adherence determinedness. In view of
Theorem \ref{thm:deptofpoints}, we will privilege depth for analogues
of $\prap$ and $\psap$, and denote by $\mathsf{C^{\prap}}$ and
$\mathsf{C^{\psap}}$ the full subcategories of $\mathsf{C}^{\Cap}$
whose objects are deep and fc-deep respectively. In view of Theorem\ref{thm:PimagePsapandPrap},
Theorem \ref{thm:deptofpoints}, and Lemma \ref{lem:ptcomptible}:
\begin{cor}
\label{cor:prappsapadj} The adjunction $\xymatrix{\Cap\ar@<1ex>[r]^{\mathbb{P}} & \mathsf{(C^{\Cap})^{op}}\ar@<1ex>[l]^{\pt L}_{\bot}}
$ restricts to $\xymatrix{\prap\ar@<1ex>[r]^{\mathbb{P}} & \mathsf{\mathsf{(C^{\prap})^{op}}}\ar@<1ex>[l]^{\pt L}_{\bot}}
$ and to $\xymatrix{\psap\ar@<1ex>[r]^{\mathbb{P}} & \mathsf{\mathsf{(C^{\psap})^{op}}}\ar@<1ex>[l]^{\pt L}_{\bot}}
$ whenever $\mathsf{C}$ is fully admissible.
\end{cor}

On the other hand, we will denote by $\mathsf{C^{\ap}}$ the full
subcategory of $\mathsf{Adh}_{\mathbb{F}_{0}}$ whose objects satisfy
(\ref{eq:pfAP}) and $\mathsf{C^{App}}$ the full subcategory of $\mathsf{C^{\prap}}$
whose objects satisfy (\ref{eq:pfAP}). Here, as the defining property
(\ref{eq:pfAP}) is in terms of the operator $\adh_{\Lambda}:L\to L^{V}$,
we privilege adherence-determinedness, that is, $\mathsf{C^{\ap}}$.
Indeed, $\mathsf{C^{\ap}}$-objects are determined by their operator
$\adh_{\Lambda}:L\to L^{V}$, while $\mathsf{C^{App}}$-objects do
not need to be. Proposition \ref{prop:adhtodeep} (\ref{enu:Boolean})
yields that $\mathsf{Adh}_{\mathbb{F}}$ embeds into $\mathsf{C^{\psap}}$
for $\mathsf{C}$ a category of bounded distributive lattices, such
as the categories of frames and of coframes.

Similarly, if $\mathsf{C}$ is a category of frames, then $\mathsf{Adh}_{\mathbb{F}_{0}}$
embeds into $\mathsf{C^{\prap}}$ and $\mathsf{C^{\ap}}$ into $\mathsf{C^{App}}$,
by Proposition \ref{prop:adhtodeep} (\ref{enu:frame}). Hence the
dotted arrows below are only valid when $\mathsf{C}$ is a category
of frames, while the solid ones are valid in any fully admissible
category $\mathsf{C}$ of distributive lattices. Note that asking
a category of frames to be fully admissible means that $\mathsf{C}$-morphisms
preserve arbitrary joins and meets, so that $\mathsf{C}$ is a subcategory
of $\mathsf{CH}$. As a result,
\begin{cor}
The adjunction $\xymatrix{\Cap\ar@<1ex>[r]^{\mathbb{P}} & \mathsf{(CH^{\Cap})^{op}}\ar@<1ex>[l]^{\pt L}_{\bot}}
$ restricts to $\xymatrix{\prap\ar@<1ex>[r]^{\mathbb{P}} & \mathsf{\mathsf{(\mathsf{Adh}_{\mathbb{F}_{0}})^{op}}}\ar@<1ex>[l]^{\pt L}_{\bot}}
$ and $\xymatrix{\Cap\ar@<1ex>[r]^{\mathbb{P}} & \mathsf{(C^{\Cap})^{op}}\ar@<1ex>[l]^{\pt L}_{\bot}}
$ restricts to $\xymatrix{\psap\ar@<1ex>[r]^{\mathbb{P}} & \mathsf{\mathsf{(\mathsf{Adh}_{\mathbb{F}})^{op}}}\ar@<1ex>[l]^{\pt L}_{\bot}}
$ whenever $\mathsf{C}$ is a fully admissible category of distributive
lattices.
\[
\xymatrix{ &  & \mathsf{C^{\Cap}}\\
\mathsf{Adh}_{\mathbb{F}_{0}}\ar@{^{(}->}[r]\ar@{^{(}.>}[rrd] & \mathsf{Adh}_{\mathbb{F}}\ar@{^{(}->}[r] & \mathsf{C^{\psap}}\ar@{^{(}->}[u]\\
\mathsf{C^{\ap}}\ar@{^{(}->}[u]\ar@{^{(}.>}[r] & \mathsf{C^{App}}\ar@{^{(}->}[r] & \mathsf{C^{\prap}}\ar@{^{(}->}[u]\\
\ap\ar@{^{(}->}[r]\ar@{^{(}->}[u]\ar@{^{(}->}[ur] & \prap\ar@{^{(}->}[rr]\ar@{^{(}->}[ru]\ar@{^{(}->}[luu] &  & \psap\ar@{^{(}->}[r]\ar@{^{(}->}[lluu] & \Cap\ar@{^{(}->}[uuull]
}
\]
\end{cor}

\subsection{Adherence-determinedness of spaces of points}

Let $\mathsf{WS}$ denote the full subcategory of $\mathsf{C^{Cap}}$
whose objects are weakly spatial, i.e., satisfy (\ref{eq:WS}).

We will say that a class $\mathbb{D}$ of filters (determined by its
instantiations $\mathbb{D}(L)\subset\mathbb{F}(L)$) is $\pt$-\emph{compatible
}if $\H_{\bullet}\in\mathbb{D}(\mathbb{P}\pt L)$ whenever $\H\in\mathbb{D}(L)$.
It is clear that:
\begin{lem}
\label{lem:ptcomptibleclassoffilters} The classes $\mathbb{F}$ of
all filters and $\mathbb{F}_{0}$ of principal filters are both $\pt$-compatible.
\end{lem}

\begin{thm}
\label{thm:ptofAdh} If $\mathsf{C}$ is fully admissible and $\mathbb{D}$
is $\pt$-compatible, then the space of points of a weakly spatial
$\mathsf{Adh_{\mathbb{D}}}$-object is a $\mathbb{D}$-adherence determined
convergence approach space. In particular, the adjunction $\xymatrix{\Cap\ar@<1ex>[r]^{\mathbb{P}} & \mathsf{(C^{\Cap})^{op}}\ar@<1ex>[l]^{\pt L}_{\bot}}
$ restricts to \\ $\xymatrix{\prap\ar@<1ex>[r]^{\mathbb{P}}_{\bot} & \mathsf{\mathsf{(Adh_{\mathbb{F}_{0}}\cap WS)^{op}}}\ar@<1ex>[l]^{\pt L}}
$.
\end{thm}

\begin{proof}
We want to show that if $(L,\Lambda)$ is $\mathbb{D}$-adherence
determined then 
\[
\bigwedge_{\mathbb{D}(\mathbb{P}(\pt L))\ni\D\#\F}\adh_{\pt L}\D(x)\leq\lambda_{\pt L}\F(x)=\bigvee\{v\in V:\Lambda\F^{\circ}(v)\in x\}.
\]
Because of (\ref{eq:regbulletonpts}), it is enough to show it under
the assumption that $\F=\F^{\natural_{\bullet}}$. 

Let $v\leq\bigwedge_{\mathbb{D}\ni\D\#\F}\adh_{\pt L}\D(x)$. Hence,
$v\leq\adh_{\pt L}\D(x)$ for every $\D\in\mathbb{D}(\mathbb{P}(\pt L))$
with $\D\#\F$, so that $\adh_{\Lambda}\D^{\circ}(v)\in x$ by Proposition
\ref{prop:adhdual-1}, and $\bigwedge_{\mathbb{D}\ni\D\#\F}\adh_{\Lambda}\D^{\circ}(v)\in x$
because $x$ is principal. Let $\H\in\mathbb{D}(L)$ with $\H\#\F^{\circ}$.
By Lemma \ref{lem:weaklyspatialgrill}, $\H_{\bullet}\#\F$ because
$\F=\F^{\natural_{\bullet}}$ and $L$ is weakly spatial. By Lemma
\ref{lem:grillstolift}, $(\H_{\bullet})^{\circ}\#\F^{\circ}$ and
$\H\subset(\H_{\bullet})^{\circ}$ so that $\adh_{\Lambda}(\H_{\bullet})^{\circ}\leq\adh_{\Lambda}\H$.
Moreover, $\D_{\H}=\H_{\bullet}\in\mathbb{D}(\mathbb{P}\pt L)$ because
$\mathbb{D}$ is $\pt$-compatible, so that 
\[
\bigwedge_{\mathbb{D}(\mathbb{P}\pt L)\ni\D\#\F}\adh_{\Lambda}\D^{\circ}(v)\leq\bigwedge_{\mathbb{D}(L)\ni\H\#\F^{\circ}}\adh_{\Lambda}\H(v)=\Lambda\F^{\circ}(v).
\]
Hence $x\in\Lambda\F^{\circ}(v)$, so that $v\leq\lambda_{\pt L}\F(x)$
and the conclusion follows.

The particulars follow from Lemma \ref{lem:ptcomptibleclassoffilters}.
\end{proof}

Note that in view of Proposition \ref{prop:adhtodeep} and Corollary
\ref{cor:prappsapadj}, weak spatiality is actually superfluous in
Theorem \ref{thm:ptofAdh} for $\mathbb{D}=\mathbb{F}$ when $\mathsf{C}$-objects
are distributive lattices, and for $\mathbb{D}=\mathbb{F}_{0}$ when
$\mathsf{C}$-objects are frames, like when $\mathsf{C=CH}$. However,
for $\mathbb{D}=\mathbb{F}_{0}$ in the category $\mathsf{C}$ of
coframes, weak spatiality is essential:
\begin{example}[A (non weakly spatial) $\mathsf{C^{\ap}}$-object whose space of points
fails to be a pre-approach space]
\label{exa:psapcounter}Let $\mathsf{C}$ denote the category of
coframes, which is fully admissible. Let $L=\mathcal{C}_{\mathbb{R}}$
be the coframe of closed subsets of $\mathbb{R}$. Let $c:\mathcal{C}_{\mathbb{R}}\to\mathcal{C}_{\mathbb{R}}$
be given by $c(F)=\cl_{\mathbb{R}}(F\setminus\{0\})\cup\{1\}$ if
$F\not\subset\{0\}$ and $c(F)=\emptyset$ if $F\subset\{0\}$. Let
\[
\Lambda\H(v)=\bigcap_{F\in\H^{\#}}c(F)
\]
for $v>0$ and $\Lambda\H(0)=\mathbb{R}$. Then, for every $v>0$,
$\Lambda(\uparrow_{L}\{0\})(v)\subset c(\{0\})=\emptyset$, while
$\Lambda(\uparrow_{L}\{x\})(v)=\bigcap_{x\in F}c(F)=\{x,1\}$ if $x\neq0$.
Thus $\pt L=\mathbb{R}\setminus\{0\}$ and $\{0\}^{\bullet}=\emptyset$
so that $L$ is not weakly spatial. For every $F_{0}\in L=\mathcal{C}_{\mathbb{R}}$
and $v>0$, 
\[
\adh_{\Lambda}F_{0}(v)=\bigvee_{F_{0}\in\H^{\#}}\Lambda\H(v)\leq c(F_{0}).
\]
If $F_{0}\subset\{0\}$ then $c(F_{0})=\adh_{\Lambda}F_{0}(v)=\emptyset$
. On the other hand, if $F_{0}\not\subset\{0\}$, there is $x\neq0$
with $x\in F_{0}$ and thus 
\[
\adh_{\Lambda}F_{0}(v)\geq\bigvee_{x\in F_{0}\setminus\{0\}}\Lambda(\uparrow_{L}\{x\})(v)=\bigvee_{x\in F_{0}\setminus\{0\}}\{x,1\}=\cl_{\mathbb{R}}(F_{0}\setminus\{0\})\cup\{1\}=c(F_{0})
\]
and thus $\adh_{\Lambda}F_{0}(v)=c(F_{0})$. Thus $\Lambda$ is $\mathbb{F}_{0}$-adherence
determined and satisfies (\ref{eq:pfAP}) because $c$ is idempotent,
so that, $(L,\Lambda)$ is a $\mathsf{C^{Ap}}$-object. On $\pt L=\mathbb{R}\setminus\{0\}$,
consider the filter $\F=\{(0,1)\}^{\uparrow}=\bigcap_{x\in(0,1)}\{x\}^{\uparrow}$.
For every $x\in(0,1)$, $(\{x\}^{\uparrow})^{\circ}=\uparrow_{L}\{x\}$,
hence $1\in\Lambda(\{x\}^{\uparrow})^{\circ}(1)$ and thus $\lambda_{\pt L}(\{x\}^{\uparrow})(1)=1$,
so that $\adh_{\pt L}B(1)=1$ whenever $B\in\F^{\#}$. Thus, $\bigwedge_{B\in\F^{\#}}\adh_{\pt L}B(1)=1$,
while $\lambda_{\pt L}\F(1)=0$ because $\F^{\circ}=\uparrow_{L}[0,1]$
and $\{0\}\in(\F^{\circ})^{\#}$ so that $\Lambda\F^{\circ}(1)\subset c(\{0\})=\emptyset$.
As a result, $\pt L$ fails to be a pre-approach space.
\end{example}

\subsection{(\ref{eq:pfAP}) and spaces of points}

Recalling Section \ref{subsec:Complete-Heyting-algebras} and Remark
\ref{rem:bullettopologiesonpts}, we see that the structure of $\pt L$
becomes more rigid in $\mathsf{CH^{\Cap}}$:
\begin{lem}
\label{lem:ptsinCH} Let $L$ be a $\mathsf{CH^{\Cap}}$-object and
let $\ell\in L$, $x\in ptL$, $A\subset\pt L$ and $a:=\bigwedge(\{A\}^{\uparrow})^{\circ}$.
Then
\begin{enumerate}
\item 
\[
\pt L\setminus\ell^{\bullet}=(\ell^{*})^{\bullet};
\]
\item 
\[
(\bigwedge x)^{\bullet}=\{x\};
\]
\item $a=\bigvee_{x\in A}\bigwedge x$;
\item $a^{\bullet}=A$.
\end{enumerate}
\end{lem}

\begin{proof}
(1). Because in $\mathsf{CH^{\Cap}}$ points of an object $L$ are
$\mathsf{CH}$-morphism $\varphi:L\to2$, we have in particular $\varphi(\ell^{*})=(\varphi(\ell))^{*}$,
that is, 
\begin{equation}
\ell\in x\iff\ell^{*}\notin x,\label{eq:pointsultra}
\end{equation}
 when the point is seen as a filter $x$ on $L$, and thus 
\[
\pt L\setminus\ell^{\bullet}=(\ell^{*})^{\bullet}.
\]

(2). Suppose $y\in(\bigwedge x)^{\bullet}$. Then $\bigwedge y\leq\bigwedge x$.
Moreover, $(\bigwedge y)\notin x$ for otherwise $\bigwedge y\leq\bigwedge x\leq(\bigwedge y)^{*}$,
so that $\bigwedge y\in x$ by (\ref{eq:pointsultra}). Thus $\bigwedge y=\bigwedge x$
and $x=y$.

(3). For every $x\in A$ and $\ell\in(\{A\}^{\uparrow})^{\circ}$,
$\ell\in x$ and $x$ is principal so that $a\in x$, that is, $\bigwedge x\leq a$
and thus $\bigvee_{x\in A}\bigwedge x\leq a$. The reverse inequality
follows from $\bigvee_{x\in A}\bigwedge x\in(\{A\}^{\uparrow})^{\circ}$,
which in turns follows from 
\[
(\bigvee_{x\in A}\bigwedge x)^{\bullet}=\bigcup_{x\in A}(\bigwedge x)^{\bullet}=A
\]
by Lemma \ref{lem:bulletsandcircles} and (2) above. 

(4). 
\[
a^{\bullet}=(\bigvee_{x\in A}\bigwedge x)^{\bullet}=\bigcup_{x\in A}(\bigwedge x)^{\bullet}=A.
\]
\end{proof}

Recall that by Proposition \ref{prop:adhtodeep}(\ref{enu:frame}),
$\mathsf{CH^{\ap}}$ is a subcategory of $\mathsf{CH^{App}}$.
\begin{thm}
Let $L$ be a weakly spatial $\mathsf{CH^{App}}$-object. Then $\pt L$
is an approach space. In other words, the adjunction $\xymatrix{\Cap\ar@<1ex>[r]^{\mathbb{P}} & \mathsf{(CH^{\Cap})^{op}}\ar@<1ex>[l]^{\pt L}_{\bot}}
$ restricts to \\
$\xymatrix{\ap\ar@<1ex>[r]^{\mathbb{P}}_{\bot} & \mathsf{\mathsf{(CH^{\ap p}\cap WS)^{op}}}\ar@<1ex>[l]^{\pt L}}
$.
\end{thm}

\begin{proof}
By Corollary \ref{cor:prappsapadj}, $\pt L$ is a pre-approach space.
We now want to show that $\pt L$ satisfies (\ref{eq:lambdaAP}) if
$L$ satisfies (\ref{eq:pfAP}). Under these assumptions, given $A\subset\pt L$
and $\epsilon\in V$, Corollary \ref{cor:adhAepsilon} yields (\ref{eq:diagonalifbulletisset}),
which combined with Lemma \ref{lem:ptsinCH} (4), yields (\ref{eq:lambdaAP}).
\end{proof}

We cannot relax the assumption on $\mathsf{C}$ to fully admissible.
\begin{example}[A weakly spatial $\mathsf{C^{App}}$-object, which is not a $\mathsf{C^{Ap}}$-object,
with a space of points failing to be approach, for $\mathsf{C}$ the
category of coframes ]
 Like in Example \ref{exa:psapcounter}, let $\mathsf{C}$ denote
the category of coframes and let $L=\mathcal{C}_{\mathbb{R}}$ be
the coframe of closed subsets of $\mathbb{R}$. Let $\N_{L}(x)=\uparrow_{L}\{x\}$
if $x\in\mathbb{R}\setminus\{0,1\}$, $\N_{L}(0)=\{F\in\mathcal{C}_{\mathbb{R}}:\exists\delta>0,[0,\delta]\subset F\}$,
and $\N_{L}(1)=\uparrow_{L}\{0,1\}$. Now, let $\Lambda\F(v)=\{x\in\mathbb{R}:\N_{L}(x)\subset\F\}$
if $v>0$ and $\Lambda\F(0)=\mathbb{R}$. Since $\F$ contains at
most one filter $\uparrow_{L}\{x\}$, $\Lambda\F(v)$ has at most
3 elements, hence is finite and thus belongs to $L=\mathcal{C}_{\mathbb{R}}$.
Note that $\pt L=\mathbb{R}$ and $F^{\bullet}=F$ for every $F\in L$,
so that $(L,\Lambda)$ is weakly spatial. Moreover, $(L,\Lambda)$
is deep. Indeed, for any $\mathbb{D}\subset\mathbb{F}L$, $\N_{L}(x)\subset\bigcap_{\D\in\mathbb{D}}\D$
if and only if $\N_{L}(x)\subset\D$ for every $\D\in\mathbb{D}$
so that $\Lambda(\bigcap_{\D\in\mathbb{D}}\D)(v)=\bigwedge_{\D\in\mathbb{D}}\Lambda\D(v)$
for every $v>0$. It is easily checked that 
\[
\adh_{\Lambda}F(v)=\begin{cases}
F\cup\{1\} & \text{if }0\in F\\
F & \text{if }0\notin F
\end{cases},
\]
for every $F\in\mathcal{C}_{\mathbb{R}}$ and $v>0$. In particular,
$0\in\adh_{\Lambda}F(\epsilon)$ if and only if $0\in F$ and thus
\[
\adh_{\Lambda}(\adh_{\Lambda}F(\epsilon))(v)=\adh_{\Lambda}F(v\otimes\epsilon)
\]
for every $F\in\mathcal{C}_{\mathbb{R}}$ and $v,\epsilon\in V$.
Hence $(L,\Lambda)$ is a weakly spatial $\mathsf{C^{App}}$-object.
Note, however, that it is not $\mathbb{F}_{0}$-adherence determined,
hence not a $\mathsf{C^{Ap}}$-object. Indeed, the filter $\F=\N_{\mathbb{R}}(0)\cap\mathcal{C}_{\mathbb{R}}$
satisfies $\Lambda\F(1)=\emptyset$ while $\bigcap_{C\in\F^{\#}}\adh_{\Lambda}C(1)=\{0,1\}$. 

To see that $\pt L$ is not an approach space, let $A=(0,1)\subset\pt L=\mathbb{R}$
and $\epsilon=1$ . The filter $\G\in\mathbb{FP}(\mathbb{R})$ generated
by $\{(0,\delta):0<\delta<1\}$ meshes with $A$ and $\G^{\circ}\supset\N_{L}(0)$
so that $0\in A^{(1)}$. Note that a proper filter finer than $\{0,1\}^{\uparrow}=\N_{L}(1)$
cannot mesh with $A$ so that $1\notin A^{(1)}=[0,1)$. On the other
hand, $\adh_{\pt L}A^{(1)}(1)=1$ because $0\in A^{(1)}$ and $(\{0\}^{\uparrow})^{\circ}=\uparrow_{L}\{0\}\supset\N_{L}(1)$.
Hence 
\[
\adh_{\pt L}A(1)=0<1=\adh_{\pt L}A^{(1)}(1)\otimes1
\]
and (\ref{eq:lambdaAP}) fails.
\end{example}

\section{\label{sec:Measure-of-compactness}Measure of compactness}

To illustrate how, in contrast to $\mathsf{C^{Conv}}$, index analysis
in the spirit of \cite{indextheory} can be conducted in $\mathsf{C^{Cap}}$,
let us consider a measure of compactness in $\mathsf{CH^{Cap}}$,
in which we can get results generalizing the classical results on
compactness. We include an illustrative example generalizing Proposition
\ref{prop:contimageofcompact} and deeper results such as generalizations
of Tychonoff's and Mrowka-Kuratowski theorems will be the subject
of future research.

\cite{pointfreecomp} introduced a notion of compactness in $\mathsf{C^{Conv}}$
and provided among other things a generalization of Proposition \ref{prop:contimageofcompact}
to that context. On the other hand, \cite[Theorem 11]{mynardmeasureCAP}
provides a generalization of Proposition \ref{prop:contimageofcompact}
to $\mathsf{Cap}$ in terms of measure of compactness. In particular,
if $(X,\lambda)$ is a convergence approach space, $\mathbb{D}$ is
a transferable class of (set-theoretic) filters, $\F\in\mathbb{FP}X$
and $A\subset X$ then the \emph{measure of $\mathbb{D}$-compactness
of $\F$ at $A$ }is given by 
\[
\mathsf{c}^{A}_{\mathbb{D}}(\F)=\bigwedge_{\mathbb{D}\ni\D\#\F}\bigvee_{a\in A}\adh_{\lambda}\D(a),
\]
and \cite[Theorem 11]{mynardmeasureCAP} implies 
\begin{prop}
\label{prop:measurecompundercont} If $(X,\lambda_{X})$ and $(Y,\lambda_{Y})$
are convergence approach spaces, $\mathbb{D}$ is a transferable class
of filters, $\F\in\mathbb{FP}X$, $A\subset X$ and $f\in\cont(X,Y)$
then
\[
\mathsf{c}^{f(A)}_{\mathbb{D}}(f[\F])\geq\mathsf{c}^{A}_{\mathbb{D}}(\F).
\]
\end{prop}

We are now in a position to propose a generalization of Proposition
\ref{prop:measurecompundercont} to $\mathsf{C^{Cap}}$, in a way
that extends Proposition \ref{prop:measurecompundercont} even in
$\mathsf{Cap}$.

Consider the \emph{measure of $\mathbb{D}$-compactness of $\F$ at
$\ell$ }given by
\[
\mathsf{c}^{\ell}_{\mathbb{D}}(\F)=\bigvee\{v\in V:\forall\D\in\mathbb{D}L,\D\#\F\then\ell\in(\adh\D(v))^{\#}\}.
\]

Let us call a class $\mathbb{D}$ of filters \emph{lattice-transferable
}if 
\[
\D\in\mathbb{D}L\then\varphi(\D)\in\mathbb{D}L'
\]
whenever $\varphi\in\mathsf{CH}(L,L')$. Recall the definition (\ref{eq:mphiCCap})
of $\mathsf{M}(\varphi)$.
\begin{thm}
\label{thm:imageofcompact} Let $\varphi\in\mathsf{CH}(L,L')$, $\F\in\mathbb{F}L'$
and $m'\in L'$ and $\mathbb{D}$ a lattice-transferable class of
filters. Then
\[
\mathsf{c}^{m'}_{\mathbb{D}}(\F)\otimes\mathsf{M}(\varphi)\leq\mathsf{c}^{\varphi_{\ell}(m')}_{\mathbb{D}}(\varphi_{\ell}(\F)).
\]
\end{thm}

\begin{proof}
Let $v\in V$ with $\adh_{L'}\L(v)\#m'$ whenever $\L\in\mathbb{D}L'$
with $\L\#\F$. Let now $\D\in\mathbb{D}L$ with $\D\#\varphi_{\ell}(\F)$,
equivalently, $\varphi(\D)\#\F$ by Lemma \ref{lem:Heygrill}. As
$\mathbb{D}$ is transferable, $\varphi(\D)\in\mathbb{\mathbb{D}}L'$
so that $\adh_{L'}\varphi(\D)(v)\#m'$. Moreover, in view of Theorem
\ref{thm:preserveadh}, 
\[
\adh_{L'}\varphi(\D)(v)\leq\varphi(\adh_{L}\varphi_{\ell}(\varphi(\D))(v\otimes\mathsf{M}(\varphi))\leq\varphi(\adh_{L}\D(v\otimes\mathsf{M}(\varphi)))
\]
because $\D\subset\varphi_{\ell}(\varphi(\D))$. Hence $\varphi(\adh_{L}\D(v\otimes\mathsf{M}(\varphi)))\#m'$,
equivalently, $\adh_{L}\D(v\otimes\mathsf{M}(\varphi))\#\varphi_{\ell}(m')$
which yields the conclusion.
\end{proof}

\begin{rem}
Recall that the classical case corresponds to $\varphi=f^{-1}:\mathbb{P}Y\to\mathbb{P}X$,
in which case $\varphi_{\ell}=f:\mathbb{P}X\to\mathbb{P}Y$.
\end{rem}

Translated in $\mathsf{Cap}$, we obtain the following generalization
of Proposition \ref{prop:measurecompundercont} that accounts for
maps that are not morphisms:
\begin{cor}
If $(X,\lambda_{X})$ and $(Y,\lambda_{Y})$ are convergence approach
spaces, $\mathbb{D}$ is a transferable class of filters, $\F\in\mathbb{FP}X$,
$A\subset X$ and $f:X\to Y$ then
\[
\mathsf{c}^{A}_{\mathbb{D}}(\F)\otimes\mathsf{m}(f)\leq\mathsf{c}^{f(A)}_{\mathbb{D}}(f[\F]),
\]
where $\mathsf{m}(f)$ is defined by \emph{(\ref{eq:mfCAP})}.
\end{cor}

\bibliographystyle{plain}

\end{document}